\documentclass{amsart}
\usepackage{latexsym, amssymb, enumerate, amsmath}

\usepackage{color}

\usepackage{mathtools}

\usepackage[backref=page]{hyperref} 
\renewcommand*{\backref}[1]{}
\renewcommand*{\backrefalt}[4]{%
     \ifcase #1 (Not cited.)%
     \or        (Cited on page~#2.)%
     \else      (Cited on pages~#2.)%
     \fi}

\usepackage[letterpaper,  margin=1in]{geometry}

\usepackage{amsthm}
\newtheoremstyle{myremarkstyle} 
{\topsep}                    
{\topsep}                    
{\normalfont}                   
{}                           
 {\bfseries}                   
{.}                          
{.5em}                       
{}  
\theoremstyle{myremarkstyle}

\def\R{\mathbb{R}}
\def\N{\mathbb{N}}

\def\Rn{\mathbb{R}^n}

\def\beginquote{``}%
\def\endquote{"}%

\def\Snm1{S^{n-1}}
\def\setMeasurableFct{\mathcal{A}}

\newcommand{\interior}[1]{\mbox{int } #1}
\newcommand{\dom}[1]{\mbox{dom } #1}

\newcommand{\prox}[1]{\left(I + {#1} \right)^{-1}}
\newcommand{\shrink}{\mbox{shrink}}

\newcommand{\characfct}[1]{\mathcal{I}_{#1}}

\begin{document}

\title{Algorithms for Overcoming the Curse of Dimensionality for
  Certain Hamilton-Jacobi Equations Arising in Control Theory and
  Elsewhere}

\author[J. Darbon]
{J\'er\^ome Darbon}
\address[J.\,Darbon]{CNRS / CMLA-Ecole Normale Sup\'erieure de Cachan, France.}
\email[J. Darbon]{jerome.darbon@cmla.ens-cachan.fr}

\author[S. Osher]
{Stanley Osher}
\address[S.\,Osher]{Department of Mathematics, UCLA, Los Angeles, CA 90095-1555, USA.}
\email[S. Osher]{sjo@math.ucla.edu}

\date{September 9, 2015 (Revised on March 10, 2016),
  Paper accepted in {\it Research in the Mathematical Sciences}. 
Research supported by ONR grants N000141410683, N000141210838 and DOE grant DE-SC00183838. }
\maketitle

\begin{abstract}
It is well known that time dependent Hamilton-Jacobi-Isaacs partial
differential equations (HJ PDE), play an important role in analyzing
continuous dynamic games and control theory problems. An important
tool for such problems when they involve geometric motion is the level
set method, \cite{OSe}.  This was first used for reachability problems
in \cite{MBT},\cite{MT}.  The cost of these algorithms, and, in fact,
all PDE numerical approximations is exponential in the space dimension
and time.

In \cite{darbon.13.cam}, some connections between
HJ-PDE and convex optimization in many dimensions are presented.
In this work we propose and test methods for
solving a large class of the HJ PDE relevant to optimal control
problems without the use of grids or numerical approximations. Rather
we use the classical Hopf formulas for solving initial value problems
for HJ PDE \cite{Ho}.  We have noticed that if the Hamiltonian is
convex and positively homogeneous of degree one (which the latter is for all
geometrically based level set motion and control and differential game
problems) that very
fast methods exist to solve the resulting optimization problem.  This
is very much related to fast methods for solving problems in
compressive sensing, based on $\ell_1$ optimization
\cite{GO},\cite{YOGD}. We seem to obtain methods which are polynomial in
the dimension.
Our algorithm is very fast, requires very low memory and is
  totally parallelizable.
We can evaluate the solution and its gradient in very high
dimensions at $10^{-4}$ to $10^{-8}$ seconds per evaluation on a
laptop.

We carefully explain how to compute numerically the optimal
  control from the numerical solution of the associated initial valued
  HJ-PDE for a class of optimal control problems. We show that our
  algorithms compute all the quantities we need to obtain easily the
  controller.

In addition, as a step often needed in this procedure, we have
developed a new and equally fast way to find, in very high dimensions,
the closest point $y$ lying in the union of a finite number of compact
convex sets $\Omega$ to any point $x$ exterior to the $\Omega$. We can
also compute the distance to these sets much faster than Dijkstra type
\beginquote fast methods \endquote, e.g. \cite{Di}.

The term \beginquote curse of dimensionality\endquote, was coined by
Richard Bellman in 1957, \cite{RB1}\cite{RB2}, when considering
problems in dynamic optimization.
\end{abstract}


\section{Introduction to Hopf Formulas, HJ PDEs and Level Set Evolutions}

We briefly introduce Hamilton-Jacobi equations with initial data and
the Hopf formulas to represent the solution. We give some examples to
show the potential of our approach, including examples to perform level
set evolutions.

Given a continuous function $H:\Rn \rightarrow \R$ bounded from below
by an affine function, we consider the HJ PDE
\begin{align}
  \label{1.1a}
\frac{\partial \varphi}{\partial t}(x,t) +
H(\nabla_x \varphi(x,t)) \;=\; 0  \quad  \text{in } \Rn \times (0,+\infty),
\end{align}
where $\frac{\partial \varphi}{\partial t} $ and $\nabla_x \varphi$
respectively denote the partial derivative with respect to $t$ and the
gradient vector with respect to $x$ of the function $\varphi$.  We are
also given some initial data
\begin{align}
\label{1.1b}
\varphi(x, 0) \;=\; J(x) & \quad  \forall\, x \in \Rn,
\end{align}
where $J:\Rn \to \R$ is convex. For the sake of simplicity we only consider
functions $\varphi$ and $J$ that are finite everywhere. Results
presented in this paper can be generalized for $H:\Rn\to \R \cup \{
+\infty\}$ and $J:\Rn \to \R \cup \{ +\infty\}$ under suitable
assumptions. We also extend our results to an interesting class
  of nonconvex initial data in section \ref{subsec.extension_future_work}.

We wish to compute the viscosity solution \cite{CL},\cite{CEL} for a
given $x \in \Rn$ and $t > 0$.

Using numerical approximations is essentially impossible for $n \geq
4$.  The complexity of a finite difference equation is exponential in
$n$ because the number of grid points is also exponential in $n$.
This has been found to be impossible, even with the use of
sophisticated e.g. ENO, WENO, DG, methods
\cite{OS},\cite{HS},\cite{ZS}. High order accuracy is no cure for
this curse of dimensionality.

We propose and test a new approach, borrowing ideas from convex
optimization, which arise in the $\ell_1$ regularization convex
optimization \cite{GO},\cite{YOGD} used in compressive sensing
\cite{CRT},\cite{Do}.  It has been shown experimentally that these
$\ell_1$ based methods converge quickly when we use Bregman and split
Bregman iterative methods. These are essentially the same as Augmented
Lagrangian methods \cite{He} and Alternating Direction Method of
Multipliers methods \cite{GM}.  These and related first order and
splitting techniques have enjoyed a renaissance since they were
recently
used very successfully for these $\ell_1$ and related problems
\cite{YOGD},\cite{GO}.
One explanation for their rapid convergence is the ``error forgetting"
property discovered and analyzed in \cite{OY} for $\ell_1$
regularization.

We will solve the initial value problem \eqref{1.1a}-\eqref{1.1b}
without discretization, using the Hopf formula~\cite{Ho}
\begin{align}
  \varphi(x,t) \;=\;& \left( J^* + t H \right)^*(x)   \label{1.2}\\
  \;=\;
  &
  -\min_{v\in R^{n}}
  \left\{ J^*(v) + tH(v) - \langle x,v \rangle\right\}
  \label{eq.Hopf_with_min}
\end{align}
where the Fenchel-Legendre transform $f^*:\Rn \to \R \cup \{
+\infty\}$
of a convex, proper, lower semicontinuous function $f:\Rn \to \R \cup
\{ +\infty\}$ is defined by
\cite{ekeland.76.book,hiriart-lemarechal.96.book-vol2,rockafellar.70.book}
\begin{equation}
  J^*(v) = \sup_{x\in R^{n}}\left\{\langle v,x \rangle - J(x)\right\},
  \label{1.3}
\end{equation}
where $\langle \cdot, \cdot\rangle$ denotes the $\ell_2(\Rn)$ inner
product. We also define for any $v \in \Rn$, $\|v\|_p =
\left(\sum_{i=1}^n |v_i|^p\right)^{\frac{1}{p}}$ for $1 \leq p < +
\infty$ and $\|v\|_\infty = \max_{i=1,\dots,n} |v_i|$.
 Note that since $J:\Rn \to \R$ is convex we have that $J^*$ is
 1-coercive \cite[Prop. 1.3.9,
 p. 46]{hiriart-lemarechal.96.book-vol2},
 i.e., $\lim_{\|x\|   \to +    \infty} \frac{J^*(x)}{\|x\|_2}\;=\;
 +\infty$.
When the gradient $\nabla_x \varphi(x,t)$ exists, then it is precisely
the unique minimizer of (\ref{eq.Hopf_with_min}); in addition,
  the gradient will also provide the optimal control considered in
  this paper (see
  Section \ref{sec.introduction_optimal_control}) For instance, this
holds for any $x\in \Rn$ and $t>0$ when $H : \Rn \to \R$ is convex,
$J:\Rn\to\R$ is strictly convex, differentiable and 1-coercive.
The Hopf formula only requires the convexity of $J$ and the continuity
of $H$, but we will often require that $H$ in (\ref{1.1a}) is also
convex.

We note that the case $H = \|\cdot\|_1$ corresponds to
$H(\nabla_x\varphi(x,t)) = \sum_{i=1}^n
\left|\frac{\partial}{\partial x_i}\varphi(x,t)\right|$, used, e.g. to
compute the Manhattan distance to the zero level set of $x\mapsto
J(x)$, \cite{Di}. This optimization is closely related to the $\ell_1$
type minimization \cite{CRT,Do},
\begin{equation}
\min_{v\in\Rn} \left\{\|v\|_1 + \frac{\lambda}{2} \|Av -
  b\|_2^2\right\},  \nonumber
\end{equation}
where $A$ is an $m \times n$ matrix with real entries, $m < n$, and
$\lambda > 0 $ arising in compressive sensing.
Let $\Omega \subset \Rn$ be a closed convex set. We denote by
$\interior \Omega$ the interior of $\Omega$.
If $J(x) < 0$ for  $ x \in \interior{\Omega}$,
 $J(x) > 0$ for $x \in (\Rn \setminus \Omega)$ and $J(x) = 0$ for
$x \in  (\Omega \setminus \interior{\Omega})$,
then the solution $\varphi$ of \eqref{1.1a}-\eqref{1.1b} has the
property that for any $t>0$ the set $\{x \in \Rn \;|\; \varphi(x,t) =
0\}$ is precisely the set of points in $\Rn \setminus \Omega$ for
which the Manhattan distance to $\Omega$ is equal to $t > 0$.  This is
an example of the use of the level set method \cite{OSe}.

Similarly, if we take $H = \|\cdot\|_2$ then for any $t>0$ the
resulting zero level set of $x \mapsto \varphi(x,t)$ will be the
points in $\Rn \setminus \Omega$ whose Euclidean distance to $\Omega$
is equal to $t$. This fact will be useful later when we find the
projection from a point $x \in \Rn$ to a compact, convex set $\Omega$.

We present here two somewhat simple but illustrative examples to show the
potential power of our approach. Time results are presented in Section
\ref{sec.numerical_results} and show that we can compute solution of
some HJ-PDEs in fairly high dimensions at a rate below a millisecond
per evaluation on a standard laptop.

Let $H = \|\cdot\|_1$ and $J(x) = \frac{1}{2} \left(\sum_{i=1}^n
  \frac{x_i^2}{a_i^2} - 1\right)$ with $a_i > 0$ for $i=1,\dots,n$.
Then, for a given $t > 0$, the set $\{x \in \Rn \;|\; \varphi(x,t) =
0\}$ will be precisely the set of points at Manhattan distance $t$
outside of the ellipsoid determined by $\{ x \in \Rn \;|\; J(x) \leq
0\}$. Following \cite[Prop. 1.3.1,
p. 42]{hiriart-lemarechal.96.book-vol2} it is easy to
see that $J^*(x) \;=\;
\frac{1}{2} \sum_{i=1}^n a_i^2 x_i + \frac{1}{2}$.
So:
\begin{align}
  \varphi(x,t) \;=\;& -\frac{1}{2} - \min_{v \in R^{n}} \left\{\frac{1}{2}
    \sum_{i=1}^n a_i^2 v_i^2 + t\sum_{i=1}^n |v_i| - \langle
    x,v\rangle\right\}
  \nonumber \\
  \;=\;& -\frac{1}{2} - \min_{v \in R^{n}} \left\{\frac{1}{2}
    \sum_{i=1}^n a_i^2 \left(v_i - \frac{x_i}{a_i^2}\right)^2 +
    t \sum_{i=1}^n |v_i|\right\} 
  + \frac{1}{2} \sum_{i=1}^n \frac{x_i^2}{a_i^2}. \nonumber
\end{align}
We note that the function to be minimized decouples into scalar
minimizations of the form
\begin{equation}
\nonumber
\min_{y\in\Rn} \left(\frac{1}{2} \|y-x\|_2^2 + \alpha\|y\|_1\right), \
\ \alpha > 0.
\end{equation}
The unique minimizer is the classical soft thresholding operator
\cite{lions.79.sna,figueiredo.98.asolimar,daubechies.04.cpam} defined
for any component~$i=1,\dots,n$ by
\begin{equation}
\label{eq.shrink1_def}
\left(\shrink_1(x,\alpha) \right)_i \;=\;
\begin{dcases}
  x_i - \alpha & \mbox{ if }  x_i > \alpha,  \\
  0         & \mbox{ if }    |x_i| \leq \alpha,  \\
  x_i + \alpha & \mbox{ if }  x_i < -\alpha.
\end{dcases}
\end{equation}
Therefore, for any $x\in\Rn$, any $t>0$ and any $i=1,\dots,n$ we have
\begin{displaymath}
\frac{\partial\varphi}{\partial x_i} (x,t)
\;=\;
\frac{1}{(a_i)^2} \ \left(\shrink_{1}(x,t)\right)_i,
\end{displaymath}
and
\begin{displaymath}
\varphi(x,t) = -\frac{1}{2} +
\sum_{i \in \{0,\dots,n\} \setminus B(t)}
\frac{1}{2}\left(\frac{|x_i|-t}{a_i}\right)^2.
\end{displaymath}
Here $B(t) \subseteq \{0,\dots,n\}$ consists of indices for which
$|x_i| \leq t$, and thus $\{0,\dots,n\} \setminus B(t)$ corresponds to
indices for which $|x_i| > t$, and the zero level set moves outwards
in this elegant fashion.

We note that in the above case we were able to compute the solution
analytically and the dimension $n$ played no significant role.  Of
course this is rather a special problem, but this gives us some idea
of what to expect in more complicated cases, discussed in
section \ref{sec.algo_optimal_control}.

We will often need another shrink operator, i.e., when we solve the
optimization problem with~$\alpha > 0$ and $x \in \Rn$
\begin{displaymath}
\min_{v \in \Rn} \left\{\frac{1}{2} \|v - x\|_2^2 + \alpha \|v\|_2\right\}.
\end{displaymath}
Its unique minimizer is given by
\begin{equation}
  \label{eq.shrink2_def}
\hbox{shrink}_2 (x,\alpha) \;=\;
\begin{cases}
 \frac{x}{\|x\|_2} \max (\|x\|_2 - \alpha, 0) & \mbox{ if } x\neq 0,
  \\
0 & \mbox{ if } x=0.
\end{cases}
\end{equation}
and thus its optimal value corresponds to the Huber function
(see \cite{winkler.03.book} for instance)
\begin{equation}
\nonumber
\min_{v \in \Rn} \left\{\frac{1}{2} \|v - x\|_2^2 + \alpha \|v\|_2\right\}
\;=\;
\begin{cases}
\frac{1}{2} \|x\|_2  & \hbox{ if } \|x\|_2 \leq \alpha, \\
 \alpha \|x\|_2 - \frac{\alpha^2}{2} & \hbox{ if } \|x\|_2 > \alpha.
\end{cases}
\end{equation}
To move the unit sphere outwards with normal velocity 1,
we use the following formula
 \begin{align}
 \varphi(x,t) &= -\min_{v \in \Rn} \left(\frac{\|v\|_2^2}{2} - t\|v\|_2 -
   \langle x,v\rangle\right) - \frac{1}{2} \nonumber \\
 &= -\min_{v\in\Rn} \left\{\frac{1}{2} \|v-x\|_2^2 + t\|v\|_2\right\} +
 \frac{1}{2} (\|x\|_2^2 - 1) \nonumber \\
 &=
 \begin{cases}
  \frac{1}{2} (\|x\|_2 - t)^2 -\frac{1}{2} & \mbox{ if }   \|x\|_2 > t\\
  -\frac{1}{2} &  \mbox{ if }  \|x\|_2 \leq t \nonumber,
 \end{cases}
 \nonumber \\
 \end{align}
 and, unsurprisingly, the zero level set of $x \mapsto \varphi(x,t)$
 is the set $x$ satisfying  $\|x\|_2 = t+1$, for $t>0$.

These two examples will be generalized below so
that we can, with extreme speed, compute the signed distance, either
Euclidean, Manhattan or various generalizations, to the boundary of
the union of a finite collection of compact convex sets.

The remainder of this paper is organized as follows:
Section~\ref{sec.introduction_optimal_control} contains
an introduction to optimal control and its connection to HJ-PDE.
Section~\ref{sec.algo_optimal_control} gives the details of our
numerical methods. Section~\ref{sec.numerical_results} presents
numerical results with some details.
We draw some concluding remarks and give future plans in
Section~\ref{sec.conclusion}. The appendix links our approach to the
concepts of gauge and support functions in convex analysis.


\section{Introduction to Optimal Control}
\label{sec.introduction_optimal_control}

First, we give a short introduction to optimal control and its
connection to HJ PDE which is given in \eqref{eq.HJ_optimal_value}.
We also introduce positively homogeneous of degree one Hamiltonians
and describe their relationship to optimal control problems. We
explain how to recover the optimal control from the solution of the
HJ-PDE. An
appendix describes further connections between these Hamiltonians and
gauge in convex analysis. Second, we present some extensions of our
work.

\subsection{Optimal control and HJ-PDE}

We are largely following the discussion in \cite{Dol}, see also
\cite{Ev}, about optimal control and its link with HJ PDE. We briefly
present it formally and we specialize it to the cases considered in
this paper.

Suppose we are given a fixed terminal time $T \in \R$, an initial time
$t < T$ along with an initial $x\in\Rn$. We consider the Lipschitz
solution $x:[t, T] \to \Rn$ of the following ordinary differential
equation
\begin{equation}
\label{eq.ode-optimal-control}
\begin{cases}
  \frac{d \mathrm{x}}{d s}(s)
  \;=\;
  f\left(\beta(s)\right)  & \mbox{ in } (t, T),\\
  \mathrm{x}(t) \;=\; x,
\end{cases}
\end{equation}
where $f: C \to \Rn$ is a given bounded Lipschitz function and $C$
some given compact set of~$\Rn$. The solution of
\eqref{eq.ode-optimal-control} is affected by the measurable function
$\beta:(-\infty, T] \to C$ which is called a control.
We note
$\setMeasurableFct = \{ \beta:(-\infty, T]\to C \; |\; \beta \mbox{ is
  measurable}\}$.
We consider the functional for given initial time $t<T$, $x\in\Rn$ and
 control $\beta$
$$
K(x,t ; \beta) \;=\;
\int_{t}^{T} L(\beta(s))ds
\,+\, J(\mathrm{x}(T)),
$$
where $\mathrm{x}$ is the solution of \eqref{eq.ode-optimal-control}.
We assume that the terminal cost $J:\Rn \to \R$ is convex. We also
assume that the running cost $L: \Rn \to \R \cup\{ + \infty\}$ is
proper, lower semicontinuous, convex, 1-coercive and $\dom{L}
\subseteq C$ where $\dom{L}$ denotes the domain of $L$ defined by
$\dom{L}\,=\, \{x\in\Rn\,|\, L(x) < +\infty \}$.  The minimization
of $K$ among all possible controls defines the value function $v:\Rn
\times (-\infty,\, T] \to \R$ given for any $x\in\Rn$ and any $t < T$ by
\begin{equation}
\label{eq.value_function}
  v(x,t) \;=\;
\inf_{\beta \in \setMeasurableFct} K(x,t; \beta).
\end{equation}
The value function \eqref{eq.value_function} satisfies the dynamic
programming principle for any $x \in \Rn$, any $t \geq T$ and any
$\tau \in (t,T)$
$$
v(x,t) \;=\;
\inf_{\beta \in \setMeasurableFct}
\left\{ \int_t^{\tau} L(\mathrm{\beta(s)})\,ds +
v\left(\mathrm{x}\left(\tau\right), \tau\right)
\right\}.
$$
The value function $v$ also satisfies the following
  Hamilton-Jacobi-Bellman equation with terminal value
\begin{equation}
\nonumber
\begin{dcases}
\frac{\partial v}{\partial t}(x,t) +
 \min_{c \in C}\{ \langle \nabla_x v(x,t), c \rangle + L(c)\} \;=\; 0
& \quad  \text{in } \Rn \times (-\infty, T),
\\
v(x,T) \;=\; J(x)  & \quad \forall x \in \Rn.
\end{dcases}
\end{equation}
Note that the control $\beta(t)$ at time $t \in  (-\infty,T)$ in
\eqref{eq.ode-optimal-control} satisfies $\langle \nabla_x v(x,t),
\beta(t) \rangle + L(\beta(t))
\;=\; \min_{c\in   C}\{\nabla_x v(x,t), c \rangle + L(c)\}$
whenever $v(\cdot,t)$ is differentiable.

Consider the function $\varphi: \Rn \times [0, +\infty) \to \R$
defined by $\varphi \left(x,t \right) \;=\; v(x,T-t)$.  We have that
$\varphi$ is the viscosity solution of
\begin{equation}
\label{eq.HJ_optimal_value}
\begin{dcases}
\frac{\partial \varphi}{\partial t}(x,t) +
H(\nabla_x \varphi(x,t)) \;=\; 0  & \quad  \text{in } \Rn \times (0,+\infty),
\\
\varphi(x,0) \;=\; J(x)  & \quad \forall x \in \Rn.
\end{dcases}
\end{equation}
where the Hamiltonian $H: \Rn \to \R \cup \{ +\infty\}$ is defined by
\begin{equation}
\label{eq.def_Hamiltonian_optimal_control}
H(p) \;=\; \max_{c \in C} \left\{ \langle -f(c), p\rangle  - L(c)\right\}.
\end{equation}
We note that the above HJ-PDE is the same as the one we consider
thoughout this paper.
In this paper we use the Hopf formula (\ref{1.2}) to solve
\eqref{eq.HJ_optimal_value}. We wish the Hamiltonian $H:\Rn\to\R$ to be not
only convex but also positively 1-homogeneous, i.e., for any $p\in \Rn$
and any $\alpha > 0$
$$
H(\alpha\, p) \;=\; \alpha \, H(p).
$$
We proceed as follows. Let us first introduce the
characteristic function
$\mathcal{I}_\Omega:\Rn\to \R\cup\{+\infty\}$ of the set $\Omega$
which is defined by
\begin{equation}
\label{eq.def_characteristic_function_set}
\characfct{\Omega}(x)  \;= \;
\begin{cases}
0 & \mbox{ if } x \in \Omega,\\
+\infty & \mbox{ otherwise}.
\end{cases}
\end{equation}
We recall that $C$ is a compact convex set of $\Rn$. We take $f(c)
\,=\, -c$ for any $c \in C$ in \eqref{eq.ode-optimal-control} and
$$
L \;=\; \mathcal{I}_C.
$$
Then, \eqref{eq.def_Hamiltonian_optimal_control}  gives the
Hamiltonian $H:\Rn\to \R$ defined by
\begin{equation}
\label{eq.H_as_support_function}
H(p) \;=\; \max_{c \,\in\, C}\, \langle c, p \rangle.
\end{equation}
Note that the right-hand side of \eqref{eq.H_as_support_function}
is called the support function of the
closed nonempty convex set~$C$ in convex analysis (see e.g.,
\cite[Def. 2.1.1, p. 208]{hiriart-lemarechal.96.book-vol1}).  We check
that $H$ defined by \eqref{eq.H_as_support_function} satisfies our
requirement. Since $C$ is bounded, we can invoke \cite[Prop. 2.1.3, p.
208]{hiriart-lemarechal.96.book-vol1} which yields that the
Hamiltonian is indeed finite everywhere. Combining \cite[Def. 1.1.1,
p. 197]{hiriart-lemarechal.96.book-vol1}
and \cite[Prop. 2.1.2, p. 208]{hiriart-lemarechal.96.book-vol1}
we obtain that $H$ is positively 1-homogeneous and convex. Of course,
the Hamiltonian can also be expressed in terms of Fenchel-Legendre
transform; we have for any $p \in \Rn$
$$
H(p) \;=\; \max_{c \,\in\, C}\, \langle c, p \rangle \;=\;
\left(\characfct{C}\right)^*(p),
$$
where we recall that the Fenchel-Legendre is defined by \eqref{1.3}.
The nonnegativity of the Hamiltonian is related to the fact that $C$
contains the origin, i.e., $0 \in C$, and gauges. This connection is
described in the appendix.

Note that the controller $\beta(t)$ for $t \in (-\infty,T)$
  in \eqref{eq.ode-optimal-control} is recovered for the solution
  $\varphi$ of \eqref{eq.HJ_optimal_value} since we have
$$
\max_{c\in C}\langle c, \nabla \varphi(x,T-t) \rangle \;=\; \langle
\beta(t), \nabla_x \varphi(x,T-t) \rangle
$$
whenever $\varphi(\cdot,t)$ is differentiable. For any $p \in \Rn$
such that $\nabla H(p)$ exists we also have $H(p) = \langle p, \nabla
H(p) \rangle$. Thus we obtain that the control is given by $\beta(t) =
\nabla H(\nabla_x \varphi(x,T-t))$.

We present in Section\ref{sec.algo_optimal_control} our
efficient algorithm that computes not only $\varphi(x,t)$ but also
$\nabla_x\varphi(x,t)$. We emphasize that we do {\it not} need to use
some numerical approximations to compute the spatial gradient. In other
words our algorithm computes all the quantities we need to get the
optimal control without using any approximations.

It is sometimes convenient to use polar coordinates. Let us denote the
$(n-1)$-sphere by $\Snm1 \;=\; \{ x \in \Rn \;|\; \|x\|_2 =1 \}$.
The set $C$ can be described in terms of the Wulff shape~\cite{OM} by
the function $W:\Snm1 \to \R$. We set
\begin{equation}
\label{2.7}
C \;=\; \{ (R\, \theta) \in \Rn \;|\;
R \geq 0,\,
\theta \in \Snm1,\,
R \leq W\left( \theta \right)\}.
\end{equation}
The Hamiltonian $H$ is then naturally defined via $\gamma: \Snm1 \to
\R$ for any $R>0$ and any $\theta \in \Snm1$ by
\begin{equation}
  \label{3.2}
H(R \theta) \; =\; R
\gamma\left(\theta \right),
\end{equation}
and where $\gamma$ is defined by
\begin{equation}
\label{eq.gamma}
\gamma(\theta)
\;=\;
\sup\left\{
\langle W(\theta')\,\theta',  \theta\rangle
\;|\;
\theta' \in \Snm1 
\right\}.
\end{equation}
The main examples are $H\;=\; \| \cdot\|_p$ for $p \in
[1,+\infty)$ and $H \;=\; \|\cdot \|_\infty$. Others include
the following two:  $H = \sqrt{\langle \cdot, A \cdot\rangle} \,=\, \|
\cdot\|_A$ with $A$ a symmetric positive definite matrix, and $H$
defined as follows for any $p \in \Rn$
$$H(p) =
\begin{cases}
 \frac{\langle p, A p  \rangle }{\|p\|_2} & \mbox{ for } p\neq 0\\
 0 & \mbox{ otherwise}.
\end{cases}
$$
In future work, we will also consider Hamiltonians defined as the
supremum of linear forms such as those that arise in linear
programming.

We will devise very fast, low memory, totally parallelizable and
apparently low time complexity  methods for solving
\eqref{eq.HJ_optimal_value} with $H$ given by
\eqref{eq.H_as_support_function} in the next section.

\subsection{Some extensions and future work }
\label{subsec.extension_future_work}

In this section we show that we can solve the problem for a much more
general class of Hamiltonians and initial data  which arise in optimal
control, including an interesting class of nonconvex initial data.

Let us first consider Hamiltonians that correspond to linear controls.
Instead of \eqref{eq.ode-optimal-control}, we consider the following
ordinary differential equation
$$
\frac{d \mathrm{x}}{d s}(s) \;=\; M\, \mathrm{x}(s) + N(s) \beta(s),
$$
where $M$ is a $n \times n$ matrices with real entries and $N(s)$ for
any $s \in (-\infty,T]$ is a  $n \times m$ matrices
with real entries. We can make a change of variables
$$
z(s) \;=\; e ^{-s M } \mathrm{x}(s),
$$
and we have
$$
\frac{d z}{d s}(s) \;=\; e ^{-s M } N(s)\beta(s).
$$
The resulting Hamiltonian now depends on $t$ and
\eqref{eq.def_Hamiltonian_optimal_control} becomes:
$$
H(p,t) \;=\;
- \inf_{c \,\in\, \Rn} \left( \langle e ^{-t M} N(t) c, p\rangle +
  L(c)\right).
$$
If $C$ is a closed convex set, and $L\; =\; \characfct{C}$ we have
$$
H(p,t) \;=\;
- \min_{c\,\in\, C}
\langle e ^{-t M } N(t) c, p\rangle,
$$
which leads to a positively 1-homogeneous convex Hamiltonians as a
function of $p$ for any fixed $t \geq 0$. If we
have convex initial data there is a simple
generalization of the Hopf formulas \cite{lions-rochet.86.ams}
\cite[Section 5.3.2, p. 215]{kurzhanski-varaiya.14.book}
$$
\varphi(z,t) \;=\;
- \min_{u \in \Rn} \left(
J^*(u) + \int_0^t H(u,s) ds - \langle z, u\rangle
\right).
$$
We intend to program and test this in our next paper.

\bigskip
We now review some well known results about the types of convex
initial value problems that yield
to max/min-plus algebra for optimal control problems (see e.g.
\cite{fleming.97.pisa,mceneaney.06.book,akian.06.book_chapter} for
instance).
Suppose we have $k$ different initial value problems $i=1,\dots,k$
\begin{equation}
\nonumber
\begin{dcases}
\frac{\partial \phi_i}{\partial t}(x,t) + H(\nabla_x \phi_i (x,t)) \;=\;
0  & \quad  \text{in } \Rn \times (0,+\infty),
\\
  \phi_i(x,0) \;=\; J_i(x)  & \quad \forall x \in \Rn. \nonumber
\end{dcases}
\end{equation}
where all initial data $J_i:\Rn\to \R$ are convex,
the Hamiltonian $H:\Rn \to R$ is convex and 1-coercive.
Then, we may use the Hopf-Lax formula to get, for any $x \in \Rn$ and
any $t > 0$
$$
\phi_i(x,t) = \min_{z \in \Rn}
\left\{
  J_i(z) + t H^*\left( \frac{x-z}{t}\right)
\right\},
$$
so
$$
\min_{i=1,\dots,k} \phi_i(x,t) =
\min_{z \in \Rn}
\left\{ \min_{i=1,\dots,k}
\left\{ J_i(z) + t H^*\left( \frac{x-z}{t}\right) \right\}
\right\}.
$$
So we can solve the initial value problem
\begin{equation}
\nonumber
\begin{dcases}
\frac{\partial \phi}{\partial t}(x,t) + H(\nabla_x \phi (x,t)) \;=\;
0  & \quad  \text{in } \Rn \times (0,+\infty),
\\
  \varphi_i(x,0) \;=\; \min_{i = 1,\dots,k} J_i(x)  & \quad \forall x \in \Rn. \nonumber
\end{dcases}
\end{equation}
by simply taking the pointwise minimum over the $k$ solutions
$\phi_i(x,t)$, each of which has convex initial data.  See
  Section~\ref{sec.numerical_results} for numerical results. As an
important example, suppose
$$
H\; =\; \| \cdot\|_2,
$$
and each $J_i$ is a level set function for a convex
compact set $\Omega_i$ with nonempty interior and where the interiors
of each $\Omega_i$ may overlap with each other.
We have $J_i(x) < 0 $ inside $\Omega_i$ , $J_i(x) >0 $ outside $\Omega_i$, and
$J_i(x) = 0$ at the boundary of $\Omega_i$. Then $\min_{i=1,\dots,n} J_i$ is
also a level set function for the union of the $\Omega_i$.
Thus we can solve complicated level set motion involving merging
fronts and compute a closest point and the associated proximal points to
nonconvex sets of this type. See section~\ref{sec.algo_optimal_control}.

\bigskip
For completeness we add the following fact about the minimum of Hamiltonians.
Let $H_i:\Rn \to \R$, with $i=1,\dots,k$, be $k$ continuous
Hamiltonians bounded from below by a common affine function.
We consider for $i=1,\dots,k$
\begin{equation}
\nonumber
\begin{dcases}
\frac{\partial \phi_i}{\partial t}(x,t) + H_i(\nabla_x \phi (x,t)) \;=\;
0  & \quad  \text{in } \Rn \times (0,+\infty),
\\
  \varphi(x,0) \;=\;  J(x)  & \quad \forall x \in \Rn, \nonumber
\end{dcases}
\end{equation}
where $J:\Rn \to \R$ is convex.
Then,
\begin{align}
  \nonumber
  \min_{i=1,\dots,k}
  \left(
    -\phi_i(x,t)
  \right)
  \;=\; &
  \nonumber
  \min_{i=1,\dots,k} \left\{
  \min_{u \in \Rn}
  \left\{
    J^*(u) + t H_i(u) - \langle u,x\rangle
  \right\}
  \right\},\\
  \;=\; &
  \nonumber
  \min_{u\in \Rn} \left\{
    J^*(u) + t \min_{i=1,\dots,k} \left\{H_i(u) - \langle u, x\rangle\right\}
  \right\},
\end{align}
that is
\begin{align}
  \max_{i=1,\dots,k} \phi_i(x,t) \;=\; -\min_{u \in \Rn}
  \left\{
    J^*(u) + t \min_{i=1,\dots,k} \left\{ H_i(u) - \langle u,x\rangle \right\}
  \right\}.
\end{align}
So we find the solution to
\begin{equation}
\nonumber
\begin{dcases}
\frac{\partial \phi}{\partial t}(x,t) + \min_{i=1,\dots,k}
H_i(\nabla_x \phi (x,t)) \;=\;
0  & \quad  \text{in } \Rn \times (0,+\infty),
\\
  \varphi(x,0) \;=\;  J(x)  & \quad \forall x \in \Rn. \nonumber
\end{dcases}
\end{equation}
by solving $k$ different initial value problems and taking the
pointwise maximum.  See Section~\ref{sec.numerical_results} for numerical results.

\bigskip
We end this section by showing that explicit formulas can be obtained
for the terminal value $\mathrm{x}(T)$ and the control $\beta(t)$ for
another class of running cost $L$.  Suppose that $C$ is a convex
compact set containing the origin and take $f(c) = -c $ for any $c \in
C$.  Assume also that $L: \Rn \to \R \cup \{+\infty\}$ is strictly
convex, differentiable when its subdifferential is nonempty and that
$\dom L$ has a non-empty interior with $\dom L \subseteq C$. Then the associated Hamiltonian $H$ is defined
by $H=L^*$.  Then, using the results of \cite{darbon.13.cam}, we have
that the $(x,t) \mapsto \varphi(x,t)$ which solves
\eqref{eq.HJ_optimal_value} is given by the Hopf-Lax formula
$\varphi(x,t) \,=\, \min_{y\in\Rn}\left\{ J(y) +
  tH^*(\frac{x-y}{t})\right\}$ where the minimizer is unique and
denoted by $\bar{y}(x,t)$. Note that the Hopf-Lax formula corresponds
to a convex optimization problem which allows us to compute
$\bar{y}(x,t)$.  In addition, we can compute the gradient with respect
to $x$ since we have $\nabla_x \varphi(x,t) = \nabla
H^*\left(\frac{x-\bar{y}(x,t)}{t}\right) \in \partial J(y(x,t))$ for any given $x\in\Rn$ and
$t>0$. For any $t \in (-\infty, T)$ and fixed $x\in\Rn$ the control is
given by $\beta(t) = \nabla H(\nabla_x \varphi(x,T-xt))$ while the
terminal value satisfies $\mathrm{x}(T) = \bar{y}(x,T-t) = x - (T-t) \nabla
H(\nabla \varphi(x,(T-t)))$. Note that both the control and the
terminal value can be easily computed. More details about these facts
will be given in a forthcoming paper.


\section{Overcoming the Curse of Dimensionality for Convex Initial
  Data and Convex Homogeneous Degree One Hamiltonians -- Optimal
  Control}
\label{sec.algo_optimal_control}

We first present our approach for evaluating the solution of the
HJ-PDE and its gradient using the Hopf formula \cite{Ho},
Moreau's identity
\cite{moreau.65.proximite} and the split Bregman algorithm \cite{GO}.
We note that the split Bregman algorithm can be replaced by other
algorithms which converge rapidly for problems of this type. An
example might be the primal-dual hybrid gradient
method \cite{zhu-chan.08.cam,chambolle.11.jmiv}. Then, we show that our
approach can be adapted to compute a closest point on a closed set,
which is the union of disjoint closed convex sets with a non empty
interior, to a given point.

\subsection{Numerical optimization algorithm}

We present the steps needed to solve
\begin{equation}
\label{3.1}
\begin{dcases}
\frac{\partial \varphi}{\partial t}(x,t) +
H(\nabla_x \varphi(x,t)) \;=\; 0  & \quad  \text{in } \Rn \times (0,+\infty),
\\
  \varphi(x,0) \;=\; J(x)  & \quad \forall x \in \Rn.
\end{dcases}
\end{equation}
We take $J:\Rn \to \R$ convex
and positively 1-homogeneous. We recall that solving \eqref{3.1},
i.e., computing
the viscosity solution, for a given $x \in \Rn,\, t > 0$ using
numerical approximations, is essentially impossible, for $n \geq 4$
due to the memory issue, and  the complexity is exponential in~$n$.

An evaluation of the solution at $x\in\Rn$ and $t>0$ for the examples
we consider in this paper is of the order of $10^{-8}$ to $10^{-4}$
seconds on a standard laptop (see Section \ref{sec.numerical_results}).
The apparent time complexity seems to be polynomial in $n$ with
remarkably small constants.

We will use the Hopf formula \cite{Ho}:
\begin{equation}
\varphi(x,t) = -\min_v \left\{J^*(v) + tH(v)-\langle
  x,v\rangle\right\}. \label{3.3}
\end{equation}
Note that the infimum is always finite and attained (i.e, it is a
minimum) since we have assumed that $J$ is finite everywhere on $\Rn$
and that $H$ is continuous and bounded from below by an affine function.

The Hopf formula \eqref{3.3} 
requires only the continuity of $H$, but we will also require
the Hamiltonian $H$ be convex as well.  We recall that the previous
section shows how to relax this
condition.

We will use the split Bregman iterative approach to solve this \cite{GO}
\begin{eqnarray}
v^{k+1} &=& \arg\min_{v \in \Rn} \{J^*(v) - \langle x,v\rangle +
\frac{\lambda}{2} \|d^k-
v-b^k\|_2^2\}, \label{3.4(a)} \\
d^{k+1} &=& \arg\min_{d \in \Rn} \left\{t H(d) + \frac{\lambda
}{2} \|d-v^{k+1}-b^k\|_2^2\right\}
\label{3.4(b)} \\
b^{k+1} &=& b^k + v^{k+1} - d^{k+1}. \label{3.4(c)}
\end{eqnarray}
For simplicity we consider $\lambda = 1$ and consider $v^0 = x$,
  $d^0=x$ and $b^0 = 0$ in this paper. The algorithm still works for
  any positive $\lambda$ and any finite values for $v^0$, $d^0$ and $b^0$.
The sequence $(v^k)_{k\in\mathbb{N}}$ and $(d^k)_{k\in\mathbb{N}}$ are
both converging to the same quantity which is a minimizer of~\eqref{3.3}.
We recall that when the minimizer of \eqref{3.3} is unique then
  it is precisely the $\nabla_x \varphi(x,t)$; in other words,
  both $(v^k)_{k\in\mathbb{N}}$ and $(v^k)_{k\in\mathbb{N}}$ converge
  to $\nabla_x \varphi(x,t)$ under this uniqueness assumption. If the
  minimizer is not unique then the sequences $(v^k)_{k\in\mathbb{N}}$
  and $(v^k)_{k\in\mathbb{N}}$ converge to an element of the
  subdifferential of $\partial (y\mapsto f(y,t))(x)$ (see below for the
  definition of a subdifferential).
We need to solve \eqref{3.4(a)} and \eqref{3.4(b)}. Note that up to
some changes of variables, both optimization problems can be
reformulated as finding the unique minimizer of
\begin{equation}
\label{eq.general_opti_pb}
\Rn \ni w \mapsto
\alpha f(w) + \frac{1}{2} \|w - z\|_2^2,
\end{equation}
where $z \in \Rn$, $\alpha > 0$, and $f:\Rn \to \R \cup \{+\infty\}$
is a convex, proper, lower semi-continuous function.  Its unique
minimizer $\bar{w}$ satisfies
the optimal condition
\begin{equation}
\nonumber
\alpha\, \partial f (\bar{w}) + \bar{w} - z \;\ni\; 0,
\end{equation}
where $\partial f (x)$ denotes the subdifferential (see
for instance \cite[p. 241]{hiriart-lemarechal.96.book-vol1},
\cite[Section 23]{rockafellar.70.book}) of $f$ at
$x \in \Rn$ and is defined by $\partial f (x) \;=\; \{ s \in \Rn\,|\,
\forall y \in \Rn,\, f(y)\, \geq\, f(x) + \langle s, y-x\rangle \}$.
We have
\begin{equation}
\nonumber
\bar{w} \;=\; \prox{\alpha\, \partial f} (z)
\;=\;
\arg\min_{w\in\Rn} \left\{ \alpha f(w) + \frac{1}{2} \|w - z\|_2^2
\right\},
\end{equation}
where $\prox{\partial f}$ denotes the \beginquote resolvent\endquote
operator of $f$ (see \cite[Def. 2, chp. 3, p. 144]{aubin.84.book},
\cite[p.  54]{brezis.73.book} for instance). It is also called the
proximal map of $f$ following the seminal paper of Moreau
\cite{moreau.65.proximite} (see also \cite[Def. 4.1.2, p.
318]{hiriart-lemarechal.96.book-vol2}, \cite[p.
339]{rockafellar.70.book}).
This mapping has been extensively studied in the context of
optimization (see for instance
\cite{combettes.07.siam-opt,eckstein-bertsekas.92.mpa,rockafellar.76.siam-co,teboulle.97.siam-op}).

\bigskip
Closed form formulas exist for the proximal of map for some specific
cases. For instance, we have seen in the introduction that
 $\prox{\alpha \ \partial \|\cdot\|_i} \;=\shrink_i(\cdot,
\alpha)$ for $i=1,2$, where we recall that $\shrink_1$  and $\shrink_2$ are
defined by \eqref{eq.shrink1_def} and \eqref{eq.shrink2_def},
respectively. Another classical example consists of considering
a quadratic form $\frac{1}{2}\|\cdot\|_A^2 \;=\; \frac{1}{2} \langle
\cdot, A \cdot \rangle$, with $A$ a symmetric positive definite
matrix with real values, which yields $\prox{\alpha\, \partial
  \left(\frac{1}{2}\|\cdot\|_A ^2\right)}\;=\; \left( I_n + \alpha \,
  A\right)^{-1}$, where $I_n$ denotes the identity matrix of size $n$.

Assume $f$ is twice differentiable with a bounded Hessian, then the
proximal map can be efficiently computed using Newton's method.
Algorithms based on Newton's method require us to solve a linear
system that involves an $n \times n$ matrix. Note that typical high
dimensions for optimal control problems are about $n=10$. For
computational purposes, these order of values for $n$ are small.

We describe an efficient algorithm to compute the proximal map of
$\|\cdot\|_\infty $ in Section \ref{subsec.prox_infty} using
parametric programming \cite[Chap. 11, Section
11.M]{rockafellar.84.book}. An algorithm to compute the proximal map
for $\frac{1}{2}\| \cdot \|_1^2$ is described in Section
\ref{subsec.l1_linfty_squared}.

\bigskip
The proximal maps for $f$ and $f^*$ satisfy the celebrated Moreau
identity \cite{moreau.65.proximite} (see also \cite[Thm. 31.5, p.
338]{rockafellar.70.book}) which reads as follows: for any $w\in \Rn$
and any $\alpha > 0$
\begin{equation}
\label{eq.moreau-identity}
\prox{\alpha\, \partial f}(w) + \alpha \prox{\frac{1}{\alpha} \, \partial f^*}
  \left( \frac{w}{\alpha}\right) \;=\; w.
\end{equation}
This shows that $\prox{\alpha \,\partial f}(w)$ can be easily computed
from $\prox{\frac{1}{\alpha} \, \partial
  f^*}\left(\frac{w}{\alpha}\right)$. In other words, depending on the
nature and properties of the mappings $f$ and $f^*$, we choose the
one for which the proximal point is \beginquote easier\endquote
to compute. Section \ref{subsec.l1_linfty_squared} describes an
algorithm to compute the proximal map of $\frac{1}{2}\|\cdot
\|_\infty^2$ using only evaluations of $\prox{ \frac{\alpha}{2} \, \partial \|\cdot\|_1^2}$
using Moreau's identity \eqref{eq.moreau-identity}.

We shall  see that Moreau's identity \eqref{eq.moreau-identity} can be
very useful to compute the proximal maps of convex and positively
1-homogeneous functions.

We consider problem \eqref{3.4(b)} that corresponds to compute the
proximal of a convex positively 1-homogeneous function $H:\Rn \to \R$.
(we use $H$ instead $f$ to emphasize that we are considering
positively 1-homogeneous functions and we set $\alpha=1$ to alleviate
notations.)  We have that $H^*$ is the characteristic function of a
closed convex set $C \subseteq \Rn$, .i.e., the Wulff shape associated
to $H$~\cite{OM},
\begin{equation}
\nonumber
H^* \;=\; \characfct{C},
\end{equation}
and $H$ corresponds to the support function $C$, that is for any $p
\in\Rn$
\begin{equation}
\nonumber
H(p)\;=\; \sup_{ s \in C}\, \langle s, p \rangle.
\end{equation}
Following Moreau \cite[Example 3.d]{moreau.65.proximite},
the proximal point of $z \in \Rn$ relative to $H^*$ is
\begin{equation}
\nonumber
\prox{\partial H^*}(z)
 \;=\;
\min_{w \in C} \left\{
  \frac{1}{2}\|w-z\|_2^2 \right\}.
\end{equation}
In other words, $\prox{ \partial H^*}(z)$ corresponds to the
projection of $z$ on the closed convex set $C$ that we denote by
$\pi_C(z)$, that is for any $z \in \Rn$
\begin{equation}
\nonumber
\prox{\partial H^*}(z)  \;=\;  \pi_C(z).
\end{equation}
Thus, using the Moreau identity \eqref{eq.moreau-identity}, we see that
$\prox{\partial H}$ can be computed from the projection on
its associated Wulff shape and we have for any $z \in \Rn$
\begin{equation}
\label{eq.prox-projection}
\prox{\partial H} (z) \;=\;
z - \pi_C(z).
\end{equation}
In other words, computing the proximal map of $H$ can be performed by
computing the projection on its associated Wulff shape
$C$. This formula is not new, see e.g. \cite{chambolle.09.ijcv,OOTT}.

Let us consider an example. Consider Hamiltonians of the form
$H = \|\cdot\|_A = \sqrt{\langle \cdot, A \cdot\rangle}$ where $A$ is a
symmetric positive matrix. Here the Wulff shape is the ellipsoid $C \;=\;
\left\{ y\in \Rn\,|\, \langle x, A ^{-1} x \rangle \; \leq 1
\right\}$. We describe in Section \ref{subsec.projection_ellipsoid} an
efficient algorithm for computing the projection on an ellipsoid.
Thus, this allows us to compute efficiently the proximal map of norms
of the form $\| \cdot \|_A$ using \eqref{eq.prox-projection}.

\subsection{Projection on closed convex set with the Level Set Method}
\label{subsec.projection}

We now describe an algorithm based on the level set method \cite{OSe}
to compute the projection $\pi_\Omega$ on a compact convex set $\Omega
\subset \Rn$ with a nonempty interior. This problem appears to be of
great interest for its own sake.

Let $\psi:\Rn \times [0,+\infty)$ be the viscosity solution of the
eikonal equation
\begin{equation}
\label{eq.eikonal}
\begin{dcases}
\frac{\partial \psi}{\partial t}(y,s) +
\|\nabla_x \psi(y,s)\|_2 \;=\; 0  & \quad  \text{in } \Rn \times (0,+\infty),
\\
  \psi(y,0) \;=\; L(y)  & \quad \forall y \in \Rn,
\end{dcases}
\end{equation}
where we recall that $\frac{\partial \psi}{\partial t}(y,t)$ and
$\nabla_x \psi(y,t)$ respectively denote the partial derivatives of
$\phi$ with respect to the time and space variable at $(y,s)$, and
where $L:\Rn \to \R$ satisfies for any $y \in \Rn$
\begin{equation}
\label{eq.L_level_set_fctx}
\begin{cases}
L(y) < 0 \mbox{ for any } y \in \mbox{int } \Omega, \\
L(y) > 0 \mbox{ for any } y \in (\Rn \setminus \Omega),\\
L(y) = 0 \mbox{ for any } y \in (\Omega \setminus \mbox{int } \Omega),\\
\end{cases}
\end{equation}
where $\mbox{int } \Omega$ denotes the interior of $\Omega$. Given
$s>0$, we consider the set
$$
 \Gamma(s) \;=\;
\left\{
y \in \Rn \,|\, \psi(y, s) \;=\;0
\right\},
$$
which corresponds to  all points that are at a (Euclidean)
distance $s$ from $\Gamma(0)$. Moreover, for a given point $y \in
\Gamma(s)$, the closest point to $y$ on $\Gamma(0)$ is exactly the
projection $\pi_\Omega(y)$ of $y$ on $\Omega$ and we have
\begin{equation}
\label{eq.projection_when_s_known}
  \pi_{\Omega}(y)
\;=\;
y  - s \frac{\nabla_x\psi(y,s)}{\|\nabla_x\psi(y,s)\|_2}.
\end{equation}
\bigskip
In this paper we will assume that $L$ is strictly convex, 1-coercive and
differentiable so that $\nabla_x \psi(s,y)$ exists for any $y \in \Rn$
and $s>0$. We note that if $\Omega$ is the finite union of sets of
this type then $\nabla_x \psi(s,y)$ may have isolated jumps. This presents
no serious difficulties.
Note that \eqref{eq.eikonal} takes the form of \eqref{3.1} with $H =
\|\cdot\|_2$ and $J = L$. We again use split Bregman to solve the
optimization given by the Hopf formula \eqref{1.2}. To avoid confusion
we respectively replace $J, v, d$ and $b$, by $L, w, e$ and $c$ in
\eqref{3.4(a)}-\eqref{3.4(c)}
\begin{eqnarray}
  w^{k+1} &=& \arg\min_{w\in\Rn} \left\{L^*(w) - \langle z,w\rangle +
    \frac{\lambda}{2} \|e^k-w -c^k\|_2^2\right\}, \label{3.9(a)}\\
  e^{k+1} &=& \arg\min_{e\in\Rn} \left\{s\|e\|_2 +
    \frac{\lambda}{2}\|e-w^{k+1}-c^k\|_2^2\right\}, \label{3.9(b)}\\
  c^{k+1} &=& c^k + w^{k+1} - e^k. \label{3.9(c)}
\end{eqnarray}
An important observation here is that $e^{k+1}$ can be solved
explicitly in \eqref{3.9(b)} using the $\shrink_2$ operator defined by
\eqref{eq.shrink2_def}.
Note that the algorithm given by \eqref{3.4(a)}-\eqref{3.4(c)} allows
us to evaluate not only $\psi(y,s)$ but
also $\nabla_x\psi(y,s)$. Indeed
for any $s>0$ $\nabla_x \psi(y,s) \;=\;
\arg\min_{v\in\Rn} \{L^*(v) + s H(v)- \langle y, v\rangle\}$, the
minimizer being unique. Thus, the above algorithm \eqref{3.9(a)}-
\eqref{3.9(c)} generates sequences $(w^k)_{k \in \mathbb{N}}$ and
$(e ^k)_{k \in \mathbb{N}}$ that both converge to $\nabla_x
\psi(z,s)$.

The above considerations about the closest point and \eqref
{eq.projection_when_s_known} give us a numerical
procedure for computing $\pi_\Omega(y)$ for any $y \in (\Rn \setminus
\Omega)$.
Find the value $\bar{s}$ so that $\psi(y, \bar{s}) \;=\; 0$ where
$\psi$ solve \eqref{eq.eikonal}. Then, compute
\begin{equation}
  \pi_{\Omega}(y)
\;=\;
y  - \bar{s} \frac{\nabla_x\psi(y,\bar{s})}{\|\nabla_x\psi(y,\bar{s})\|_2},
\end{equation}
to obtain the projection $\pi_{\Omega}(z)$. We compute $\bar{s}$ using
Newton's method to find the $0$ of the function $(0, +\infty) \ni s
\mapsto \psi(z,s)$. Given an initial $s_0 > 0$ the Newton iteration
corresponds to computing for integers $l>0$
\begin{equation}
\label{eq.Newton_zero}
s_{l+1} \;=\;
s_{l} - \psi(z,s_{l}) \left(\frac{\partial \psi}{\partial t} (z,s_{l})
\right)^{-1}.
\end{equation}
From \eqref{eq.eikonal} we have $\frac{\partial \psi}{\partial t}
(z,s) = - \|\nabla_x \psi(z,s)\|_2$ for any $s > 0$. We can thus
compute \eqref{eq.Newton_zero}.

It remains to choose the initial data $L$ related to the set
$\Omega$. We would like it to be smooth so that the proximal point in
\eqref{3.9(a)} can be computed efficiently using Newton's
method. (If $L$ lacks differentiability the approach can be
  easily modified using \cite[Chap 11, Section
  11.M]{rockafellar.84.book} or using \cite{cheng-tsai.08.jcp}.) We
consider
$\Omega$ as Wulff shapes that are expressed, thanks to \eqref{2.7},
with the function $W:\Snm1\to \R$, that is
\begin{equation}
\nonumber
\Omega \;=\; \{ (R\, \theta) \in \Rn \;|\;
R \geq 0,\,
\theta \in \Snm1,\,
R \leq W\left( \theta \right)\}.
\end{equation}
As a simple example, if $W \equiv 1$, then we might try
$L = \|y\|_2 - 1$, the signed distance to $\Gamma(0)$.  This does not
suit our purposes, because its Hessian is singular and its dual is the
indicator function of the $l_2(\Rn)$ unit ball. Instead we take
\begin{displaymath}
L(y) = \frac{1}{2} \left(\|y\|_2^2 - 1\right).
\end{displaymath}
Note $L^*(y) \; =\; \frac{1}{2} \left(\|y\|_2^2 + 1\right)$, both of
these are convex and $\mathcal{C}^2$ functions with non vanishing
gradients away from the origin, i.e. near $\Gamma(s) $.  This gives us
a hint as how to proceed.

Recall that we need to get initial data which behaves as a level set function
should; i.e., as defined by \eqref{eq.L_level_set_fctx}. We also want
either $L$ or $L^*$ to be smooth enough, actually twice differentiable
with Lipschitz continuous Hessian, so that Newton's method can be
used.

We might take
$$
L( R\, \theta)\;=\; \frac{1}{2m}
\left(
\left(\frac{R}{W\left(\theta\right)}\right)^{2m} -1
\right),
$$
where $R \geq 0$, $\theta \in \Snm1$ and $m$ a positive integer. If we
consider the important case of $W(\theta) \,=\, \|\theta\|_A
^{-\frac{1}{2}} \,=\, \frac{1}{\sqrt{\langle \theta, A\, \theta
    \rangle}}$, with $A$ a symmetric positive definite matrix, which
corresponds to $\Omega = \{x \in \Rn \,|\, \sqrt{\langle x, A\, x \rangle}
\,\leq\, 1\}$, then $m\geq2$ will have a smooth enough Hessian for $L$. In
fact, $m=2$
will lead us to a linear Hessian. There will be situations where using
$L^*$ is preferable because $L$ cannot be made smooth enough this way.

We can obtain $L^*$ using polar coordinates and taking
$1 \geq m > \frac{1}{2}$, that is for any $R\geq 0$ and any~$\theta \in \Snm1$
\begin{align}
\nonumber
L^*(R\,  \theta)
\;=\;&
\displaystyle
\sup_{r\geq 0,\, v \in \mathcal{S}^{n-1}}
\left\{
R\, r\, \langle \theta, v\rangle + \frac{1}{2m} - \frac{1}{2m}
\left( \frac{r}{W(v)} \right)^{2m}
\right\}
\\
\nonumber
\;=\;&
\frac{1}{2m} \sup_{\|v\|_2=1}
\left\{
(2m-1)
\left( R \langle \theta, v\rangle W (v) \right)^{\frac{2m}{2m-1}} +1
\right\}
\\
\nonumber
\;=\;&
\frac{1}{2m}
\left\{
\left( R\, \gamma\left( \theta \right)\right)^{\frac{2m}{2m-1}}
(2m-1) +1
\right\}
\\
\nonumber
\;=\;&
\frac{1}{2m}
\left\{
\left( H\left( R\,\theta \right) \right)^{\frac{2m}{2m-1}}
(2m-1) +1
\right\}
\end{align}
where we recall that $\gamma$ is defined by \eqref{eq.gamma}.
As we shall see, $L^*$ is often preferable to $L$ yielding a
smooth Hessian for $ 1 \geq m > \frac{1}{2}$ with $m$ close to $\frac{1}{2}$.

Let us consider some important examples. If the set $\Omega$ is defined
by $\Omega = \{x\in\Rn \,|\, \|x\|_p \leq 1\}$ for $1 < p < \infty$,
then we can consider two cases
\begin{itemize}
\item[(a)] $2 \leq p < + \infty$
\item[(b)] $1 < p \leq 2$
\end{itemize}
For case $(a)$ we need only take
$$
L \;=\; \frac{1}{2m}
\left(
\|\cdot\|_p^{2m}  -1
\right).
$$
If $m \geq 1$ it is easy to see that Hessian of $L$ is continuous and
bounded for $p\geq 2$. So we can use Newton's method for the
first choice above.
For case $(b)$, we construct the Fenchel-Legendre of the function
$$
L \;=\; \frac{1}{2m} \left( \|\cdot\|_p^{2m} - 1\right),
$$
but this time for $\frac{1}{2} < m \leq 1$. It is easy to
see that 
$$
L^* \;=\; \frac{1}{2m}
\left(
(2m-1)\|\cdot\|_q^{\frac{2m}{2m-1}} + 1
\right),
$$
where $\frac{1}{p} + \frac{1}{q} = 1$ and if $\frac{1}{2} < m \leq 1$
the Hessian of $L^*$ is continuous.

Other interesting examples include the following regions defined by
$$
\Omega \;=\; \left\{
x \in \Rn \;|\;  \langle x, Ax \rangle \leq \|x\|_2
\right\}.
$$
For $\Omega$ to be convex we require $A$ to be a positive definite
symmetric matrix with real entries and its maximal eigenvalue is
bounded by twice the minimal eigenvalue.
We can take
\begin{equation}
\label{eq.def_L_quad_over_norm}
L(x) \;= \;
\begin{cases}
\frac{1}{2m}
\left(
\left(
\frac{\langle x,\, A x
\rangle}{\|x\|_2} \right)^{2m}
-1 \right) & \mbox { if } x \neq 0,\\
0 & \mbox { if } x = 0.
\end{cases}
\end{equation}
for $m\geq 2$ and see that this has a smooth Hessian.

\section{Numerical Results}
\label{sec.numerical_results}

We will consider the following Hamiltonians
\begin{itemize}
\item $H = \|\cdot\|_p$ for $p = 1, 2, \infty$,
\item $H = \sqrt{\langle \cdot, A \cdot\rangle}$ with $A$ symmetric positive
  definite matrix,
\end{itemize}
and the following initial data
\begin{itemize}
\item $J = \frac{1}{2}\|\cdot\|_p^2$ for $p = 1, 2, \infty$,
\item $J = \frac{1}{2}\langle \cdot, A\cdot \rangle$ with $A$ a
  positive definite diagonal matrix.
\end{itemize}
It will be useful to consider the spectral decomposition of $A$, i.e.,
$A = P D P^\dagger$ where $D$ is a diagonal matrix, $P$ is an
orthogonal matrix and $P^\dagger$ denotes the transpose of $P$. The
identity matrix in $\Rn$ is denoted by $I_n$.

First, we present the algorithms to compute the proximal points for
the above Hamiltonians and initial data. We shall describe these
algorithms using the following generic formulation for the proximal map
$$
\prox{\alpha \,\partial f} (w)
 \;=\;
\arg \min_{w  \in \Rn}
\left\{
\frac{1}{2} \|w-z\|_2^2 +
\alpha f(w)
\right\}.
$$
Second, we present the time results on a standard laptop. Some time
results are also provided for a 16 cores computer which shows that our
approach scales very well. The latter is due to
our very low memory
requirement. Finally some plots that represent the
solution of some HJ PDE are presented.

\subsection{Some explicit formulas for simple specific cases}

We are able to obtain explicit formulas for the proximal map for some
specific cases of interest in this paper. For instance, as we have
seen, considering $f = \| \cdot \|_1$ gives for any $i=1,\dots,n$
$$
\left(\prox{\alpha \,\partial \|\cdot \|_1}(w)\right)_i
\;=\; \mbox{sign}(w_i) \max(|w_i|-\alpha, 0),
$$
where $\mbox{sign}(\beta) = 1$ if $\beta\geq 0$ and $-1$ otherwise.
The case $f=\|\cdot\|_2$ yields a similar formula
$$
\prox{\alpha \,\partial \|\cdot \|_2}(w)
\;=\;
\frac{\mu}{\alpha + \mu} x,
$$
with
$$
\mu = \max(\|w\|_2 - \alpha,0).
$$
The two above cases are computed in linear time with respect to the
dimension $n$.

Proximal maps for positive definite quadratic forms, i.e., $f(w) =
\frac{1}{2}\langle w, A w\rangle$ are also easy to compute since for
any $\alpha >0$ and any $z\in\Rn$
$$
\prox{ \alpha \,\partial f}(z) \;=\; ( I_n + \alpha A)^{-1} (z)
\;=\;
 P (I_n + \alpha D)^{-1} P^\dagger z.
$$
where we recall that $A = P D P^\dagger$ with $D$ a diagonal matrix
and $P$ an orthogonal matrix. The time complexity is dominated by the
evaluation of the matrix-vector product involving $P$ and
$P^{\dagger}$.
\subsection{The case of $\|\cdot \|_\infty$ }
\label{subsec.prox_infty}

Let us now consider the case $f=\|\cdot\|_\infty$. Since
$\|\cdot\|_\infty$ is a norm, its Fenchel-Legendre transform is the
indicator function of its dual norm ball,
that is $(\|\cdot\|_\infty)^* \;=\;
\mathcal{I}_C$ with $C \;=\; \{ z \in \Rn \;|\;
\|z\|_1 \leq 1\}$. We use Moreau's identity
\eqref{eq.prox-projection} to compute $\prox{\alpha \|
  \cdot\|_\infty}$; that is for any $\alpha>0$ and for any $z \in \Rn$
$$
\prox{\alpha \partial\,\| \cdot\|_\infty}(z)
\;=\; z - \alpha \pi_C(\frac{z}{\alpha}) = z - \pi_{\alpha C} (z)
$$
where we recall that $\pi_{\alpha C}$ denotes the projection operator onto
the closed convex set $\alpha C$. We use a simple variation of parametric
approaches that are well-known in graph-based  optimization algorithm
(see \cite[Chap. 11, Section 11.M]{rockafellar.84.book} for instance).

Let us assume that $z \notin (\alpha C)$. The projection corresponds
to solve
$$
\pi_{\alpha C} (z) \;=\;
\begin{dcases}
\arg \min_{w\in\Rn} \frac{1}{2}\|w-z\|_2^2 \\
\mbox{s.t. } \|w\|_1 \leq \alpha.
\end{dcases}
$$
 Now we use Lagrange duality (see
 \cite[chap. VII]{hiriart-lemarechal.96.book-vol1} for instance). The
 Lagrange dual function $g: [0,+\infty) \to \R$ is defined by
$$
     g(\mu)
     \;=\;
     \min_{w\in\Rn}\left\{ \frac{1}{2} \| z -w \|_2^2 +
       \mu \|w\|_1 - \mu \alpha \right\},
$$
that is
$$
g(\mu)
\;=\;
\| \mbox{shrink}_1(z, \mu)\|_1 - \mu \alpha.
$$
We denote by $\bar{\mu}$ the value that realizes the
maximum of $g$.  Then, we obtain for any $z
\notin \alpha C$ that
\begin{equation}
\label{eq.projection_mu_linfty}
\pi_{\alpha C}(z) \;=\;
\mbox{shrink}_1(z, \bar{\mu}),
\end{equation}
where  $\bar{\mu} $
satisfies
\begin{equation}
\label{eq.condition_mu_linfty}
\|\mbox{shrink}_1(z,\bar{\mu})\|_1 \;=\; \alpha.
\end{equation}
Computing the projection is thus reduced to computing the optimal
value of the Lagrange multiplier $\bar{\mu}$.
Consider the function $h : [0,\|z\|_1] \to \R$ defined by $ h(\mu)
\;=\; \|\mbox{shrink}_1(z, \mu)\|_1$. We have that $h$ is continuous,
piecewise affine, $h(0)=\|z\|_1$, $h(\|z\|_1) = 0$ and $h$ is
decreasing (recall that we assume $z \notin (\alpha C)$). Following
\cite[Chap. 11, Section 11.M]{rockafellar.84.book}, we call
breakpoints the values for which $h$ is
not differentiable. The set of breakpoints for $h$ is $B = \left\{0
\right\} \bigcup_{i=1,\dots,n} \left\{ |z_i|\right\}.$ We sort all
breakpoints in increasing order and we denote this sequence by $(l_i,
\dots, l_{m}) \in B^{m}$ with $l_i < l_{i+1}$ for $i=1,\dots,(m-1)$, where $m
\leq n$ is the number of breakpoints. This operation takes $O(n \log
n)$.  Then, using a bitonic search, we can find $j$ such that
$\|\shrink_1(z, l_j)\|_1 \leq \|\shrink_1(z, \bar{\mu})\|_1 <
\|\shrink_1(z, l_{j+1})\|_1$ in $O(n \log n)$. Since $h$ is affine on
$[\|\shrink_1(z, l_j)\|_1, [\|\shrink_1(z, l_{j+1})\|_1]$ a simple
interpolation computed in constant time  yields $\bar{\mu}$ that
satisfies \eqref{eq.condition_mu_linfty}. We then use
\eqref{eq.projection_mu_linfty} to compute the projection. The overall
time complexity is, therefore, $O(n \log n)$.

\subsection{The case $\|\cdot\|_A$ and projection
  on a ellipsoid}
\label{subsec.projection_ellipsoid}

We follow the same approach as for $\|\cdot\|_\infty$. We consider $f
= \| \cdot \|_A = \langle \cdot, A \cdot \rangle$ which is a norm
since $A$ is assumed to be symmetric positive definite.  The dual norm
is $\|\cdot\|$ \cite[Prop. 4.2, p.  19]{ekeland.76.book} is
$\| \cdot\|_{A ^{-1}}$.
Thus $\left(\| \cdot\|_{A}\right)^* = I_{\mathcal{E}_A}$ with
$\mathcal{E}_A$ defined by
$$
\mathcal{E}_A
\;=\; \left\{ y \in \Rn \;|\; \langle y, A ^{-1} y \rangle \leq 1\right\}.
$$
Using Moreau's identity \eqref{eq.prox-projection} we only need to
compute the
 projection $\pi_{\mathcal{E}_{A}}(w)$ of $w\in \Rn$ on the ellipsoid
$\mathcal{E}_{A} $.
Note that we have $\pi_{\mathcal{E}_{A}} (w) = P\,
\pi_{\mathcal{E}_{D }} (P^\dagger w)$ where we recall that
$A=PDP^\dagger$ with $D$ and $P$ a diagonal and orthogonal matrix,
respectively. Thus we only describe the algorithm for the projection
on an ellipsoid involving positive definite diagonal matrices.

To simplify notation we take $d_i = D_{ii}$ for $i=1,\dots,n$.
We consider the ellipsoid $\mathcal{E_D}$ defined by
$$
\mathcal{E}_{D} \;=\;
\left\{
  x \in \mathbb{R}^n\, |\, \sum_{i=1}^{n} \left(
    \frac{x_i}{d_i}\right)^2 \leq 1
 \right\}.
$$
Let $w \notin \mathcal{E}_D$.
We can easily show (see \cite[Exercise III.8]{hiriart.98.book} for
instance) that $\Pi_{\mathcal{E}_D}(w)$ satisfies for any
$i=1.\dots,n$
\begin{equation}
\label{eq.lagrange_multiplier_equation_ellipsoid}
\left(\Pi_{\mathcal{E}_D}(w)\right)_i \;=\;
\frac{d_i^2 w_i}{d_i^2 + \bar{\mu}},
\end{equation}
where the Lagrange multiplier $\bar{\mu} > 0$ is the unique solution
of $\sum_{i=1}^n \frac{d_i^2 w_i^2}{(d_i^2 + \mu)^2} = 1$. We find
such $\bar{\mu}$ by minimizing the function $[0,+\infty) \ni \mu
\mapsto \sum_{i=1}^{n} d_i ^2 w_i ^2 (d_i ^2 + \mu)^{-1} + \mu$ using
Newton's method which generates a sequence $(\mu_k)_{k \in \N}$
converging to $\bar{\mu}$. We set the initial value to $\mu_0 = 0$ and
we stop Newton's iterations for the first $k$ which statisfies
$|\mu_{k+1} - \mu_{k}| \leq 10^{-8}$. Once we have the value for $\bar{\mu}$ we
use \eqref{eq.lagrange_multiplier_equation_ellipsoid} to obtain the
approximate projection.

\subsection{The cases $\frac{1}{2}\|\cdot \|_1^2$ and
$\frac{1}{2}\|\cdot\|_\infty^2$}
\label{subsec.l1_linfty_squared}

First, we consider the case of $\frac{1}{2}\|\cdot \|_1^2$. We have
for any $\alpha >0$ and any $z\in \Rn$
$$
\prox{\frac{\alpha}{2}\, \partial \|\cdot \|_1^2} (z)
\;=\;
\prox{\alpha \|\cdot\|_1 \,\partial\|\cdot\|_1}(z).
$$
Thus, assuming there exists $\bar{\beta} \geq 0$ such that
\begin{equation}
  \label{eq.optimality_beta_l1_2}
\bar{\beta} \;=\;
\alpha \left\| \shrink_1(z,\bar{\beta}) \right\|_1,
\end{equation}
we have for any $\alpha >0$ and any $z\in \Rn$
\begin{equation}
  \label{eq.prox_l1_2_with_beta}
\prox{\frac{\alpha}{2} \,\partial \|\cdot \|_1^2} (z)
\;=\; \shrink_1(z, \bar{\beta}).
\end{equation}
The existence of $\bar{\beta}$ and an algorithm to compute
it follow.
Let us assume that $z\neq 0$ (otherwise, $\bar{\beta}=0$ works and the
solution is of course 0). Then, consider the function
$g:\left[0,\|z\|_1\right] \to \R$ defined by
$$
g(\beta)  \;=\;
\alpha \|\shrink_1(z,\beta)\|_1 - \beta.
$$
It is continuous, and $g(0)= \alpha \|z\|_1$ while
$g\left(\|z\|_1\right) = -\|z\|_1$. The intermediate value theorem
tells us that there exists $\bar{\beta}$ such that $g(\bar{\beta}) =
0$, that is, satisfying \eqref{eq.optimality_beta_l1_2}.

The function $g$ is decreasing, piecewise affine and the breakpoints
of $g$ (i.e, the points where $g$ is not differentiable) are $B= \{0\}
\cup_{i=1,\dots,n} \{ |z_i|\}$. We now proceed similarly as for the
case $\|\cdot\|_\infty$. We note $(l_i,\dots, l_{m}) \in B^{m}$ the
breakpoints sorted in increasing order, i.e., such that $l_i <
l_{i+1}$ for $i=1,\dots,(m-1)$, where $m \leq n$ is the number of
breakpoints. We use a bitonic search to find the two consecutive
breakpoints $l_i$ and $l_{i+1}$, such that
$g(l_i) \;\geq\;  0 \;>\; g(l_{i+1})$. Since $g$ is affine on $[l_i,
l_{i+1}]$ a simple interpolation yields the value $\bar{\beta}$. We
then compute $\prox{\frac{\alpha}{2} \partial \|\cdot\|_1^2}(z)$
using \eqref{eq.prox_l1_2_with_beta}.

\medskip
We now consider the case $\frac{1}{2}\|\cdot\|_\infty^2$. We have
(for instance \cite[Prop. 4.2, p.  19]{ekeland.76.book})
$$
\left(\frac{1}{2}\| \cdot \|_\infty^2\right)^*
\;=\;
\frac{1}{2}\| \cdot \|_1^2
$$
Then Moreau's identity \eqref{eq.moreau-identity} yields for any
$\alpha >0 $ and for any $z \in \Rn$
$$
\prox{\frac{\alpha}{2}\,\partial \|\cdot\|_\infty^2}(z)
\;=\;
z - \prox{\frac{\alpha}{2}\,\partial \|\cdot\|_1^2}(z),
$$
which can be easily computed using the above algorithm for evaluating
$\prox{\frac{\alpha}{2}\,\partial \|\cdot\|_1^2}(z)$.

\subsection{Time results and illustrations}

We now give numerical results for several Hamiltonians and initial
data. We present time results on a standard laptop using a
single core which show that our approach allows us to evaluate very
rapidly HJ PDE solutions. We also present some time results on a
16 cores computer to show that are approach scales very well.
We also present some plots that depict the solution of some HJ PDEs.

We recall that we consider the following Hamiltonians
\begin{itemize}
\item $H = \|\cdot\|_p$ for $p = 1, 2, \infty$,
\item $H = \sqrt{\langle \cdot, D \cdot\rangle}$ with $D$ a diagonal positive
  definite matrix,
\item $H = \sqrt{\langle \cdot, A \cdot\rangle}$ with $A$ symmetric positive
  definite matrix,
\end{itemize}
and the following initial data
\begin{itemize}
\item $J = \frac{1}{2}\|\cdot\|_p^2$ for $p = 1, 2, \infty$.
\item $J = \frac{1}{2}\langle \cdot, D^{-1}\cdot \rangle$ with $D$ a
  positive definite diagonal matrix,
\end{itemize}
where the matrix $D$ and $A$ are defined follows: $D$ is a diagonal
matrix of size $n\times n$ defined by $D_{ii} = 1 + \frac{i-1}{n-1}$
for $i=1,\dots,n$. The symmetric positive definite matrix $A$ of size
$n\times n$ is defined by $A_{ii} = 2$ for $i=1,\dots,n$
and $A_{ij} = 1 $ for $i,j=1,\dots,n$ with $i\neq j$.

All computations are performed using IEEE double precision
floating-points where denormalized number mode has been disabled.
The quanties $(x,t)$ are drawn uniformly in $[-10,10]^n \times
  [0,10]$. We present the average time to evaluate a solution for
  $1,000,000$ runs.

We set $\lambda =1$ in the split Bregman algorithm
\eqref{3.4(a)}-\eqref{3.4(c)}. We stop the iterations when the
following stopping criteria is met:
$\| v^k - v^{k-1}\|_2^2 \leq 10^{-8}$
and
$\| d^k - d^{k-1}\|_2^2 \leq 10^{-8}$ and
$\| d^k - v^{k}\|_2^2 \leq 10^{-8}$.

We first carry out the numerical experiments on an Intel Laptop Core
i5-5300U running at 2.3~GHz. The implementation here is {\it single}
threaded, i.e., only one core is used. Tables \ref{table.l2_2},
\ref{table.linfty_2}, \ref{table.l1_2} and \ref{table.diag_2} present
time results for several dimensions $n=4,8,12,16$ and with initial
data $J = \frac{1}{2}\|\cdot\|_2^2$, $J =
\frac{1}{2}\|\cdot\|_\infty^2$, $J = \frac{1}{2}\|\cdot\|_1^2$ and $J
= \frac{1}{2}\langle \cdot, D \rangle$ respectively. We see that it
takes about $10^{-8}$ to $10^{-4}$ seconds per evaluation of the
solution.

We now consider experiments that are carried out on a computer with 2
Intel Xeon E5-2690 processors running at 2.90GHz. Each processor has 8
cores. Table \ref{table.parallel} present the average time to compute
the solution with Hamiltonian $H=\|\cdot\|_\infty$ and initial data
$J=\|\cdot\|_1^2$ for
several dimensions and various number of used cores. We see that our
approach scales very well. This is due to the fact that our algorithm
requires little memory which easily fits in the L1 cache of each
processor. Therefore cores are not competing for resources. This
suggests that our approach is suitable for low-energy embedded
systems.

\begin{figure}[t]
\begin{minipage}[b]{.49\linewidth}
  \centering
 \centerline{\includegraphics[scale=0.35]{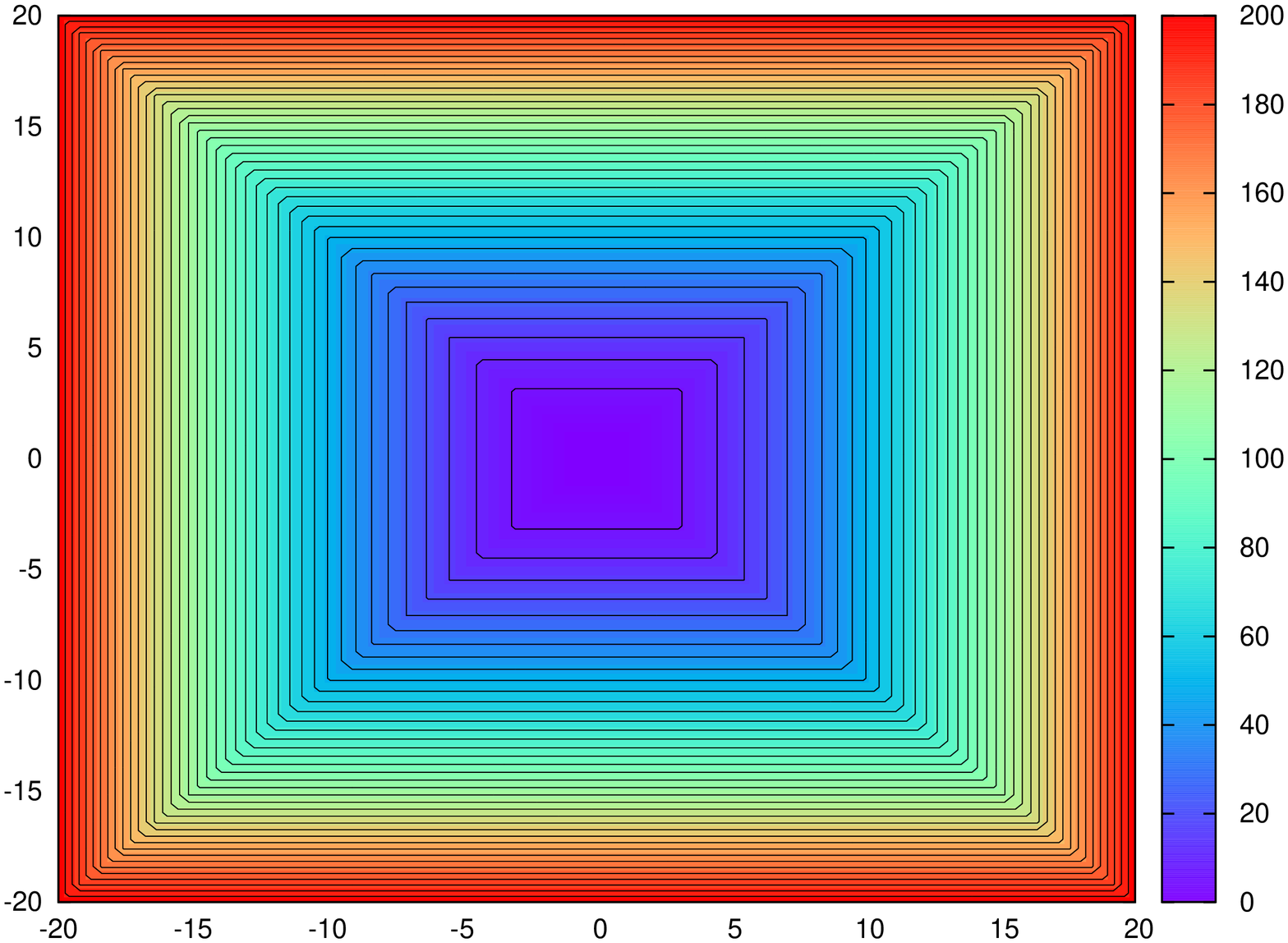}}
 \vspace{-.1cm}
  \centerline{\footnotesize{(a)}}\medskip
\end{minipage}
\hfill
 \begin{minipage}[b]{0.49\linewidth}
   \centering
  \centerline{\includegraphics[scale=0.35]{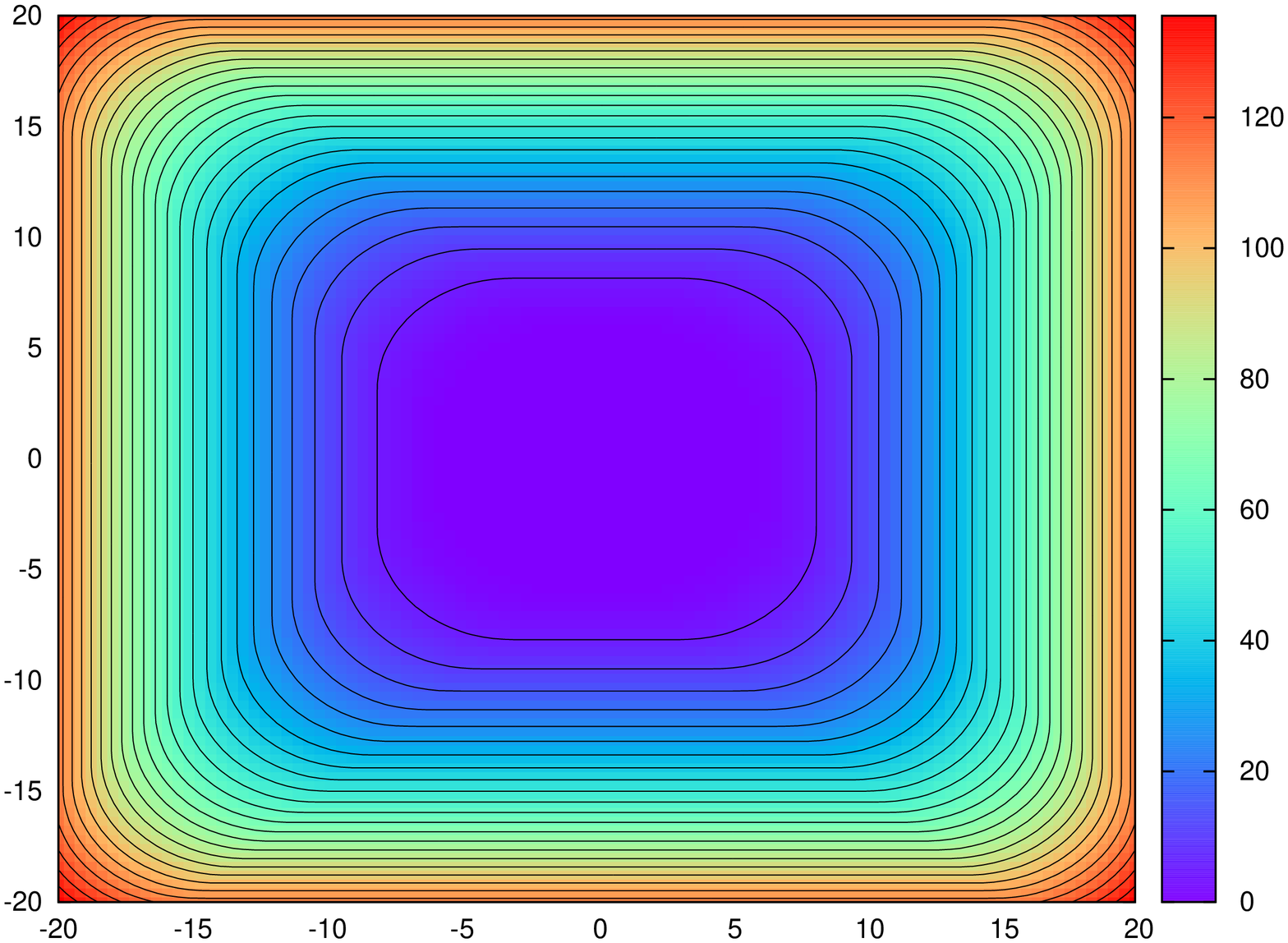}}
   \vspace{-.1cm}
   \centerline{\footnotesize{(b)}}\medskip
 \end{minipage}
\begin{minipage}[b]{0.49\linewidth}
  \centering
 \centerline{\includegraphics[scale=0.35]{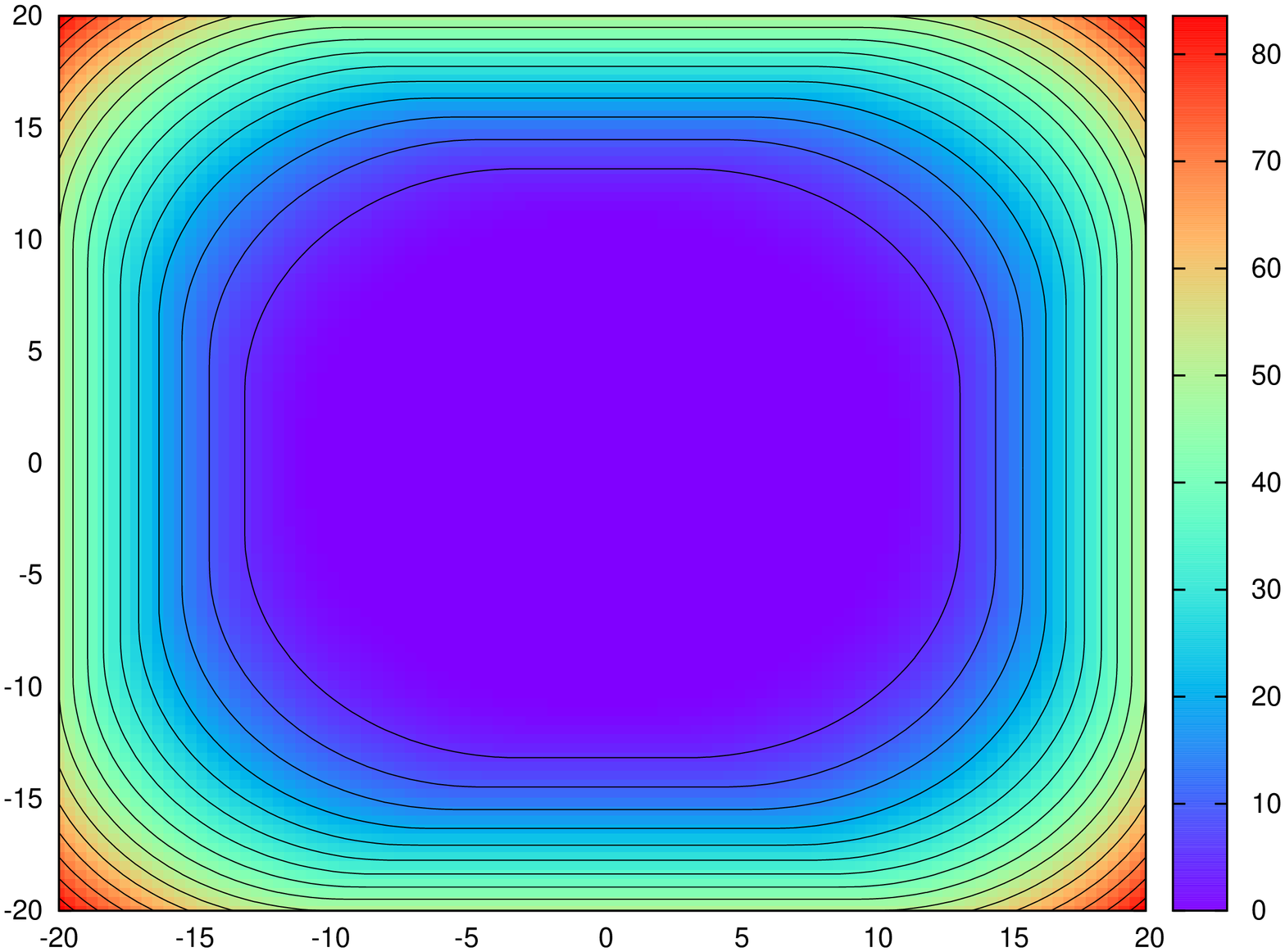}}
  \vspace{-.1cm}
  \centerline{\footnotesize{(c)}}\medskip
\end{minipage}
 \begin{minipage}[b]{0.49\linewidth}
    \centering
    \centerline{\includegraphics[scale=0.35]{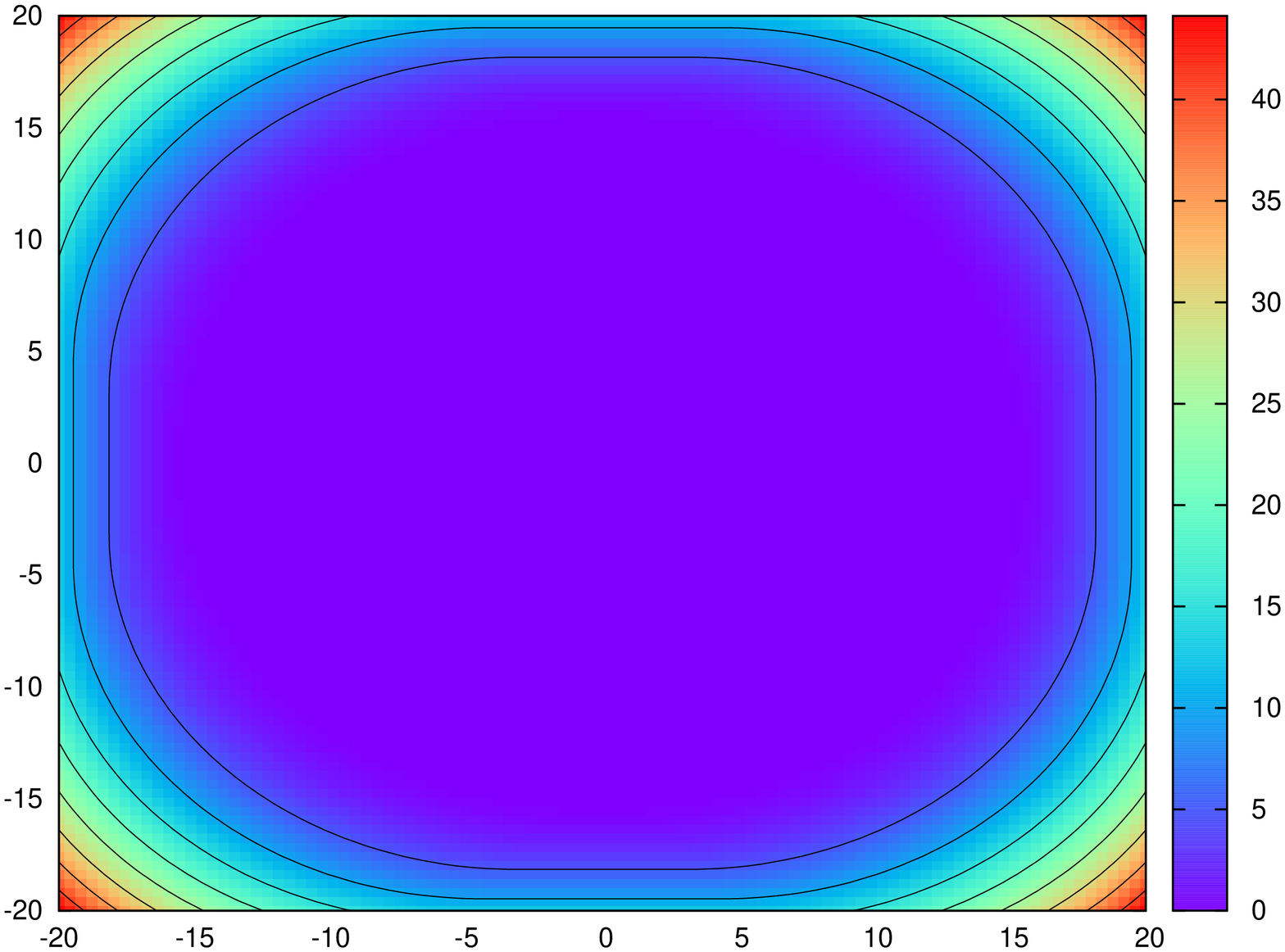}}
    \vspace{-.1cm}
    \centerline{\footnotesize{(d)}}\medskip
  \end{minipage}

  \caption{ Evaluation of the solution $\phi((x_1,x_2,
    0,0,0,0,0,0)^\dagger,t)$ of
    the HJ-PDE with initial data  $J=\frac{1}{2}\|\cdot\|_\infty^2$
    and Hamiltonian $H=\|\cdot\|_2$ for $(x_1,x_2) \,\in\, [-20,20]^2$ for
    different times $t$.
    Plots for $t=0, 5,10,15$ and respectively depicted in (a),
    (b), (c) and (d). The level lines multiple of 5 are superimposed
    on the plots.}
\label{fig.J_linfty_2__H__l2}
\end{figure}

\begin{figure}[ht]
\begin{minipage}[b]{.49\linewidth}
  \centering
 \centerline{\includegraphics[scale=0.35]{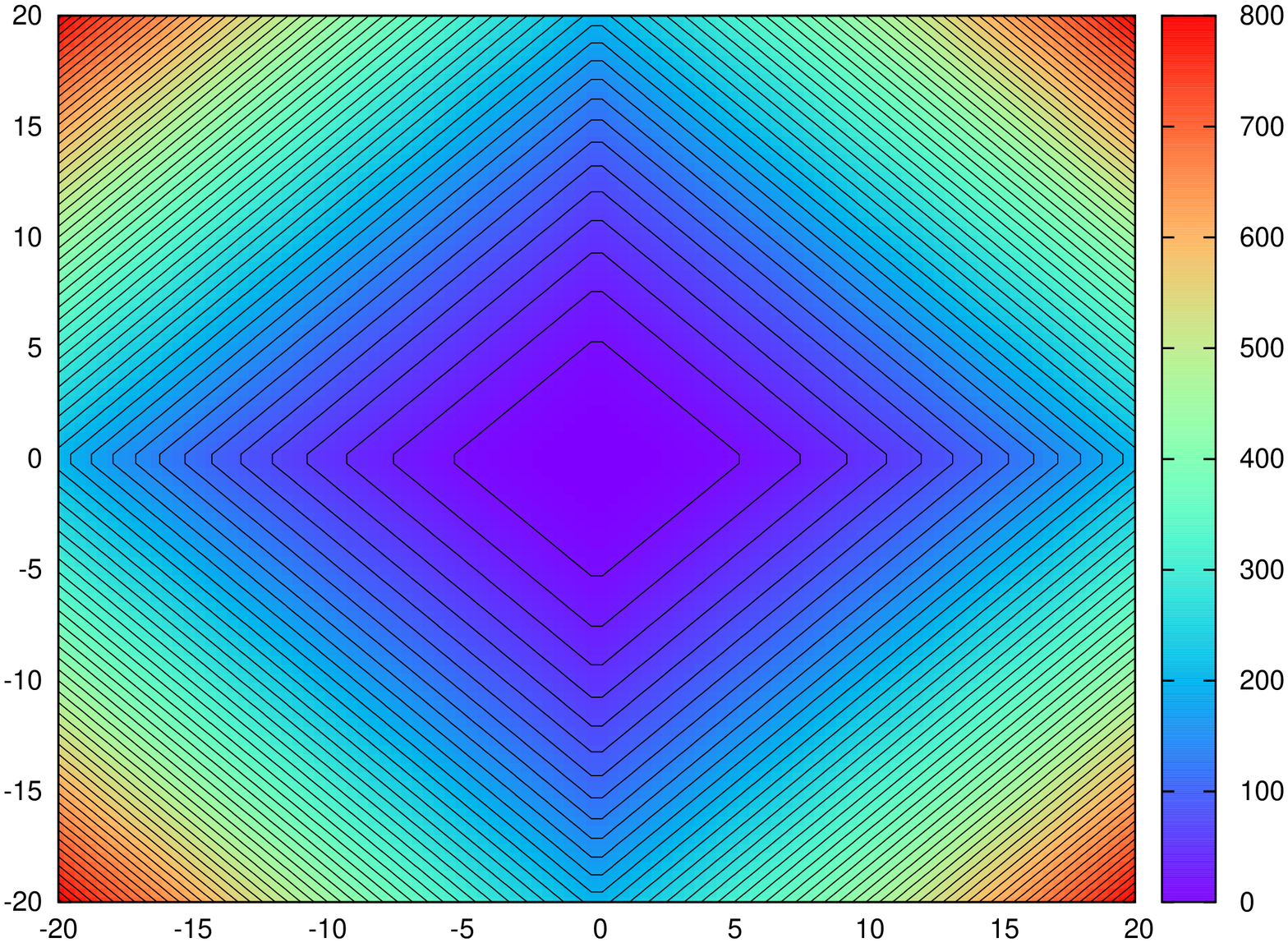}}
 \vspace{-.1cm}
  \centerline{\footnotesize{(a)}}\medskip
\end{minipage}
\hfill
 \begin{minipage}[b]{0.49\linewidth}
   \centering
  \centerline{\includegraphics[scale=0.35]{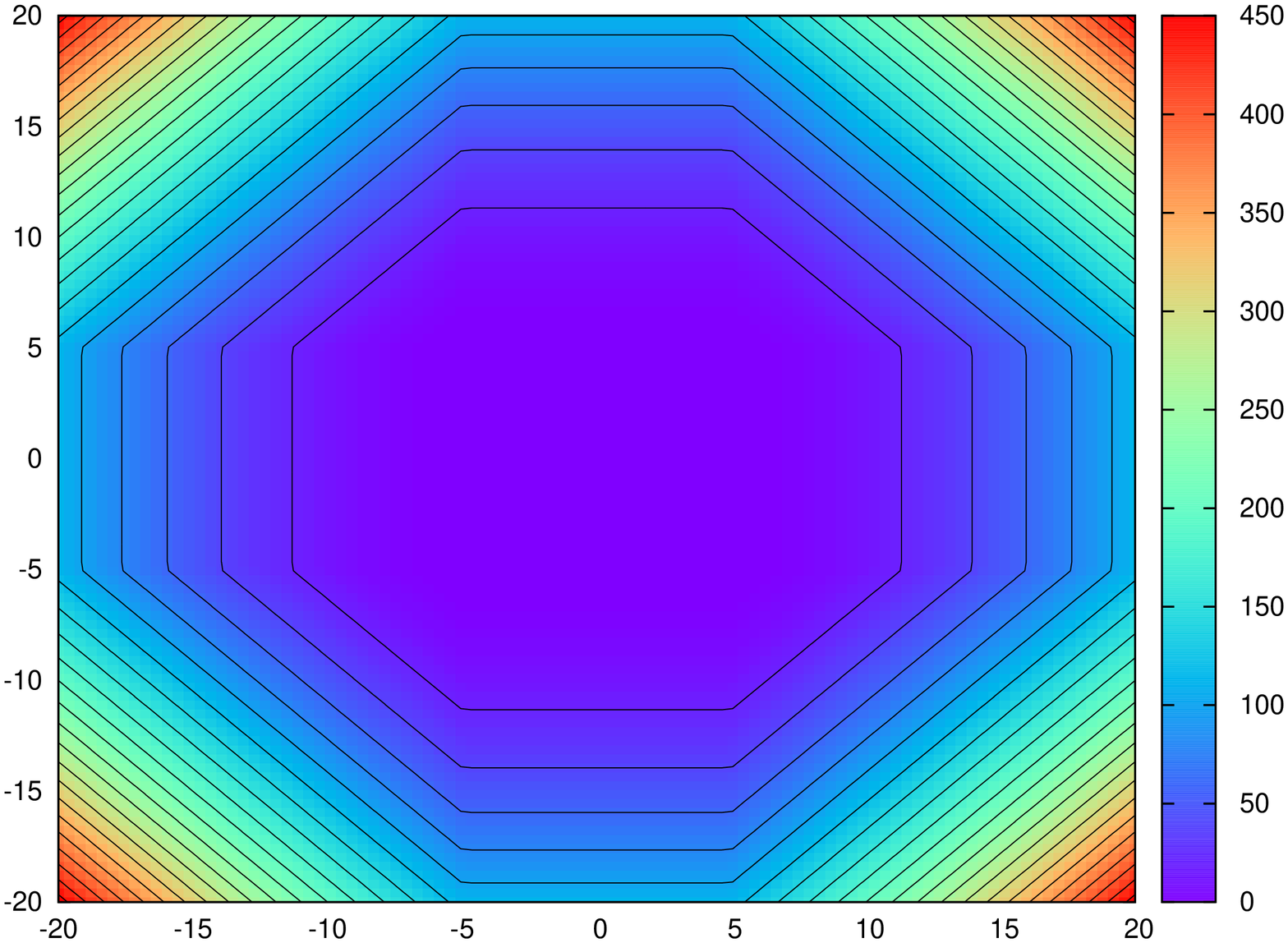}}
   \vspace{-.1cm}
   \centerline{\footnotesize{(b)}}\medskip
 \end{minipage}
\begin{minipage}[b]{0.49\linewidth}
  \centering
 \centerline{\includegraphics[scale=0.35]{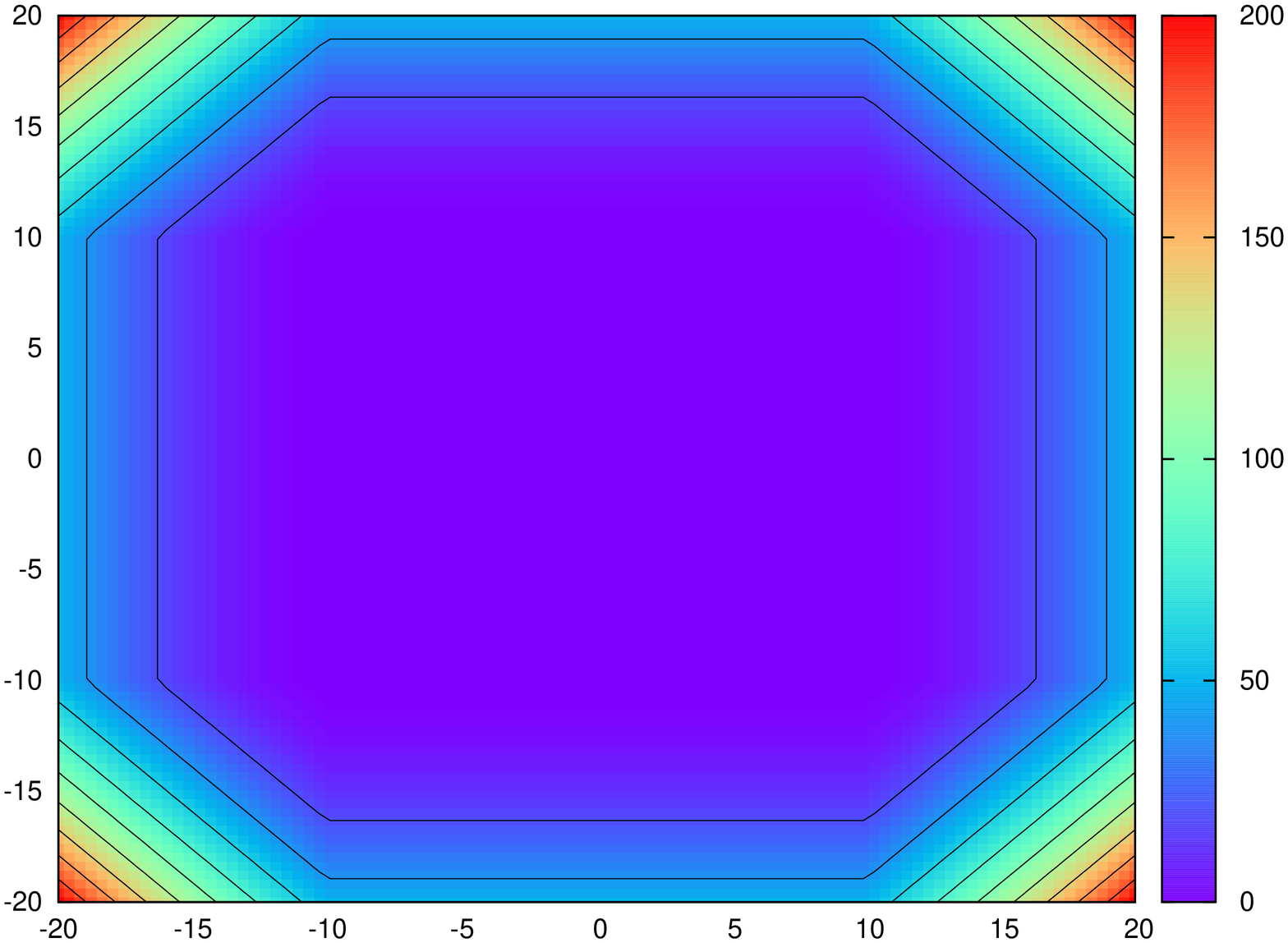}}
  \vspace{-.1cm}
  \centerline{\footnotesize{(c)}}\medskip
\end{minipage}
 \begin{minipage}[b]{0.49\linewidth}
    \centering
    \centerline{\includegraphics[scale=0.35]{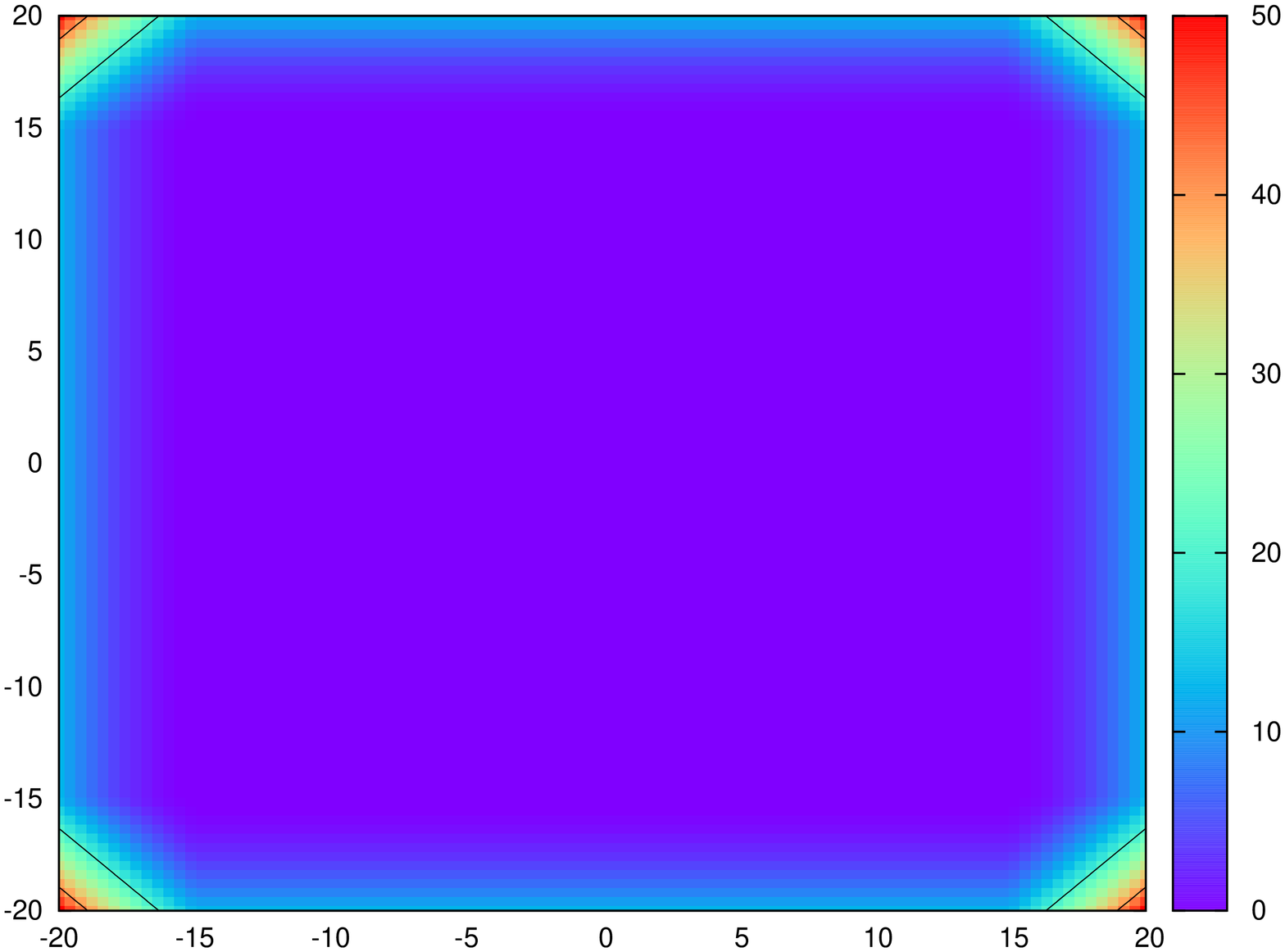}}
    \vspace{-.1cm}
    \centerline{\footnotesize{(d)}}\medskip
  \end{minipage}

  \caption{
Evaluation of the solution $\phi((x_1,x_2,0,0,0,0,0,0)^\dagger,t)$ of
    the HJ-PDE with initial data  $J=\frac{1}{2}\|\cdot\|_1^2$
    and Hamiltonian $H=\|\cdot\|_1$ for $(x_1,x_2) \,\in\, [-20,20]^2$ for
    different times $t$.  Plots for $t=0, 5,10,15$ are
    respectively depicted in (a), (b), (c) and (d). The level
    lines multiple of 20 are superimposed   on the plots.}
\label{fig.J_l1_2__H__l1}
\end{figure}

\begin{figure}[th]
\begin{minipage}[b]{.49\linewidth}
  \centering
 \centerline{\includegraphics[scale=0.35]{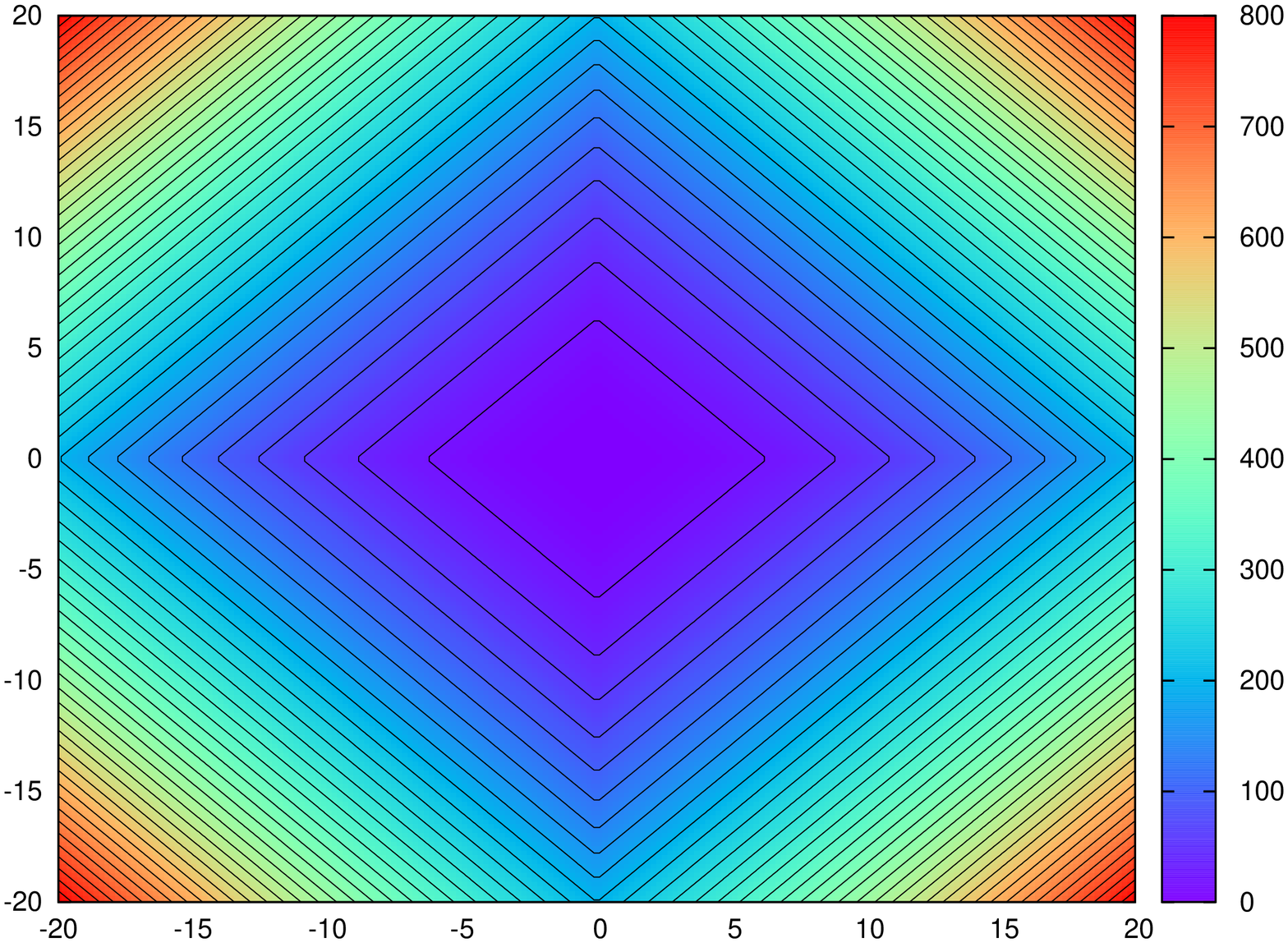}}
 \vspace{-.1cm}
  \centerline{\footnotesize{(a)}}\medskip
\end{minipage}
\hfill
 \begin{minipage}[b]{0.49\linewidth}
   \centering
   \centerline{\includegraphics[scale=0.35]{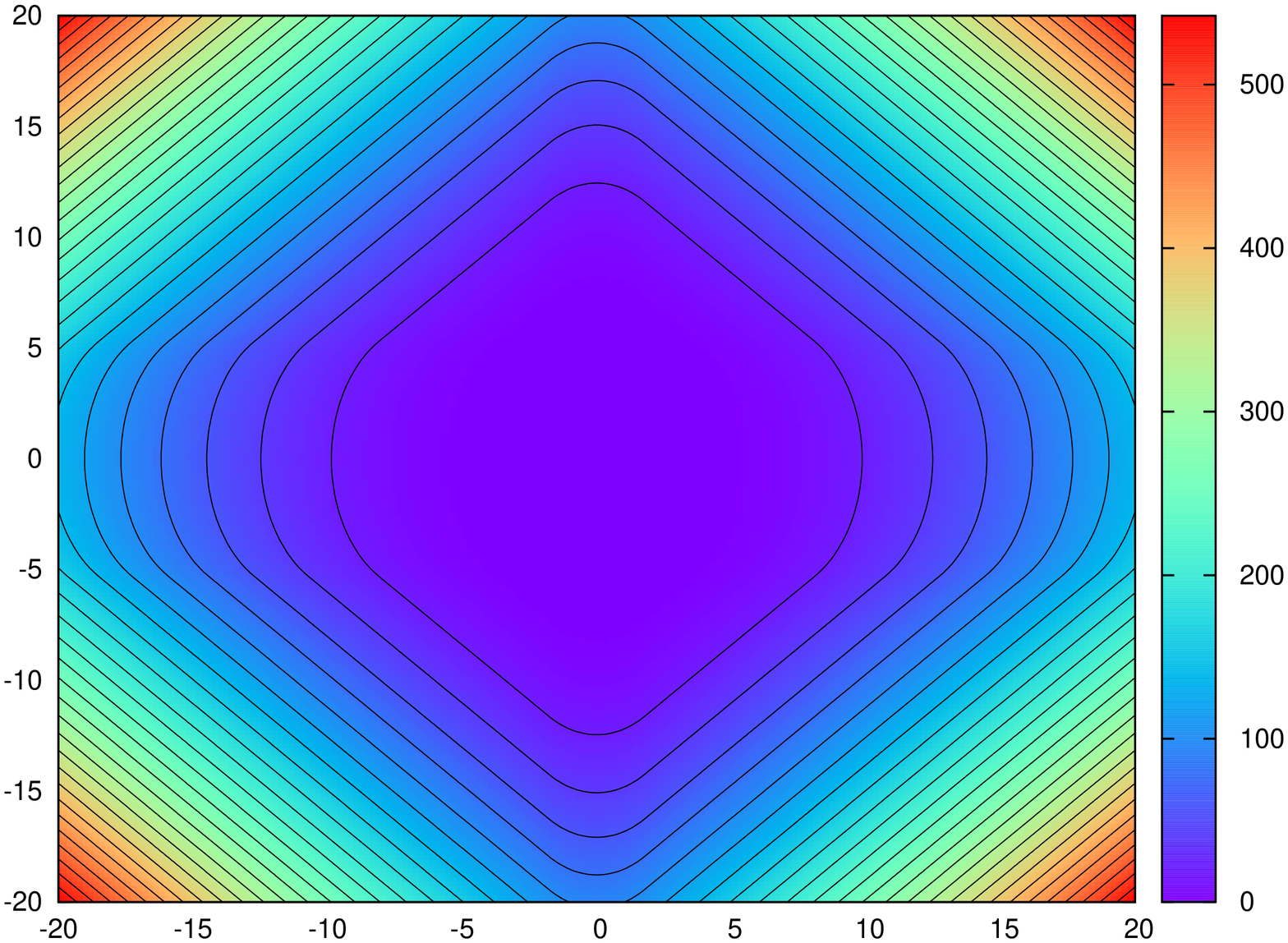}}
   \vspace{-.1cm}
   \centerline{\footnotesize{(b)}}\medskip
 \end{minipage}
\begin{minipage}[b]{0.49\linewidth}
  \centering
  \centerline{\includegraphics[scale=0.35]{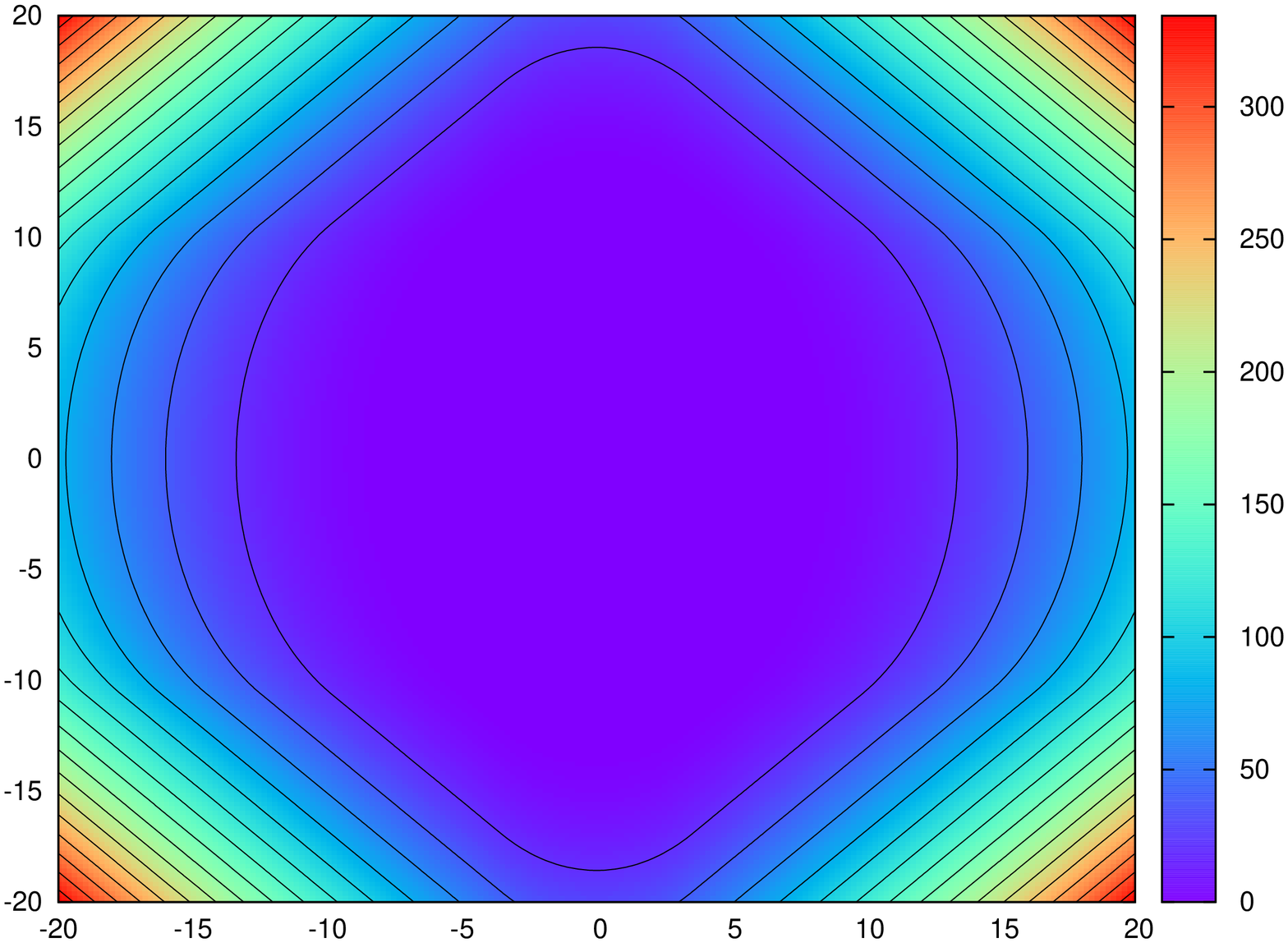}}
  \vspace{-.1cm}
  \centerline{\footnotesize{(c)}}\medskip
\end{minipage}
\begin{minipage}[b]{0.49\linewidth}
  \centering
  \centerline{\includegraphics[scale=0.35]{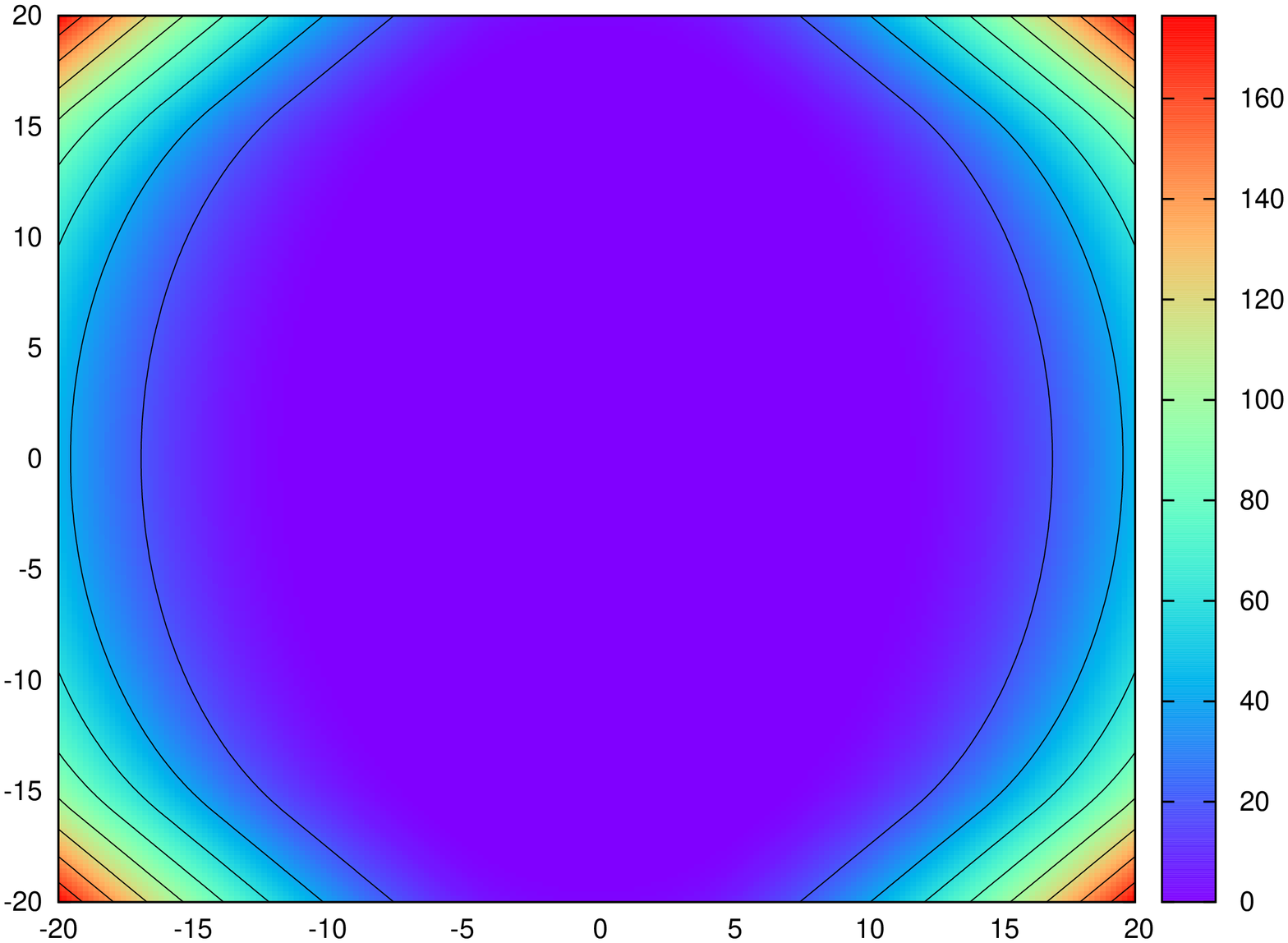}}
  \vspace{-.1cm}
  \centerline{\footnotesize{(d)}}\medskip
\end{minipage}

  \caption{
    Evaluation of the solution $\phi((x_1,x_2,0,0,0,0,0,0)^\dagger,t)$ of
    the HJ-PDE with initial data  $J=\frac{1}{2}\|\cdot\|_1^2$
    and Hamiltonian $H=\|\cdot\|_D$, for $(x_1,x_2) \,\in\,
    [-20,20]^2$ for   different times $t$.
    Plots for $t=0, 5,10,15$ are
    respectively depicted in (a), (b), (c) and (d). The level
    lines multiple of 20 are superimposed   on the plots.
}
\label{fig.J_l1_2__H__DIAG_NORM}
\end{figure}

\begin{figure}[th]
\begin{minipage}[b]{.49\linewidth}
  \centering
 \centerline{\includegraphics[scale=0.35]{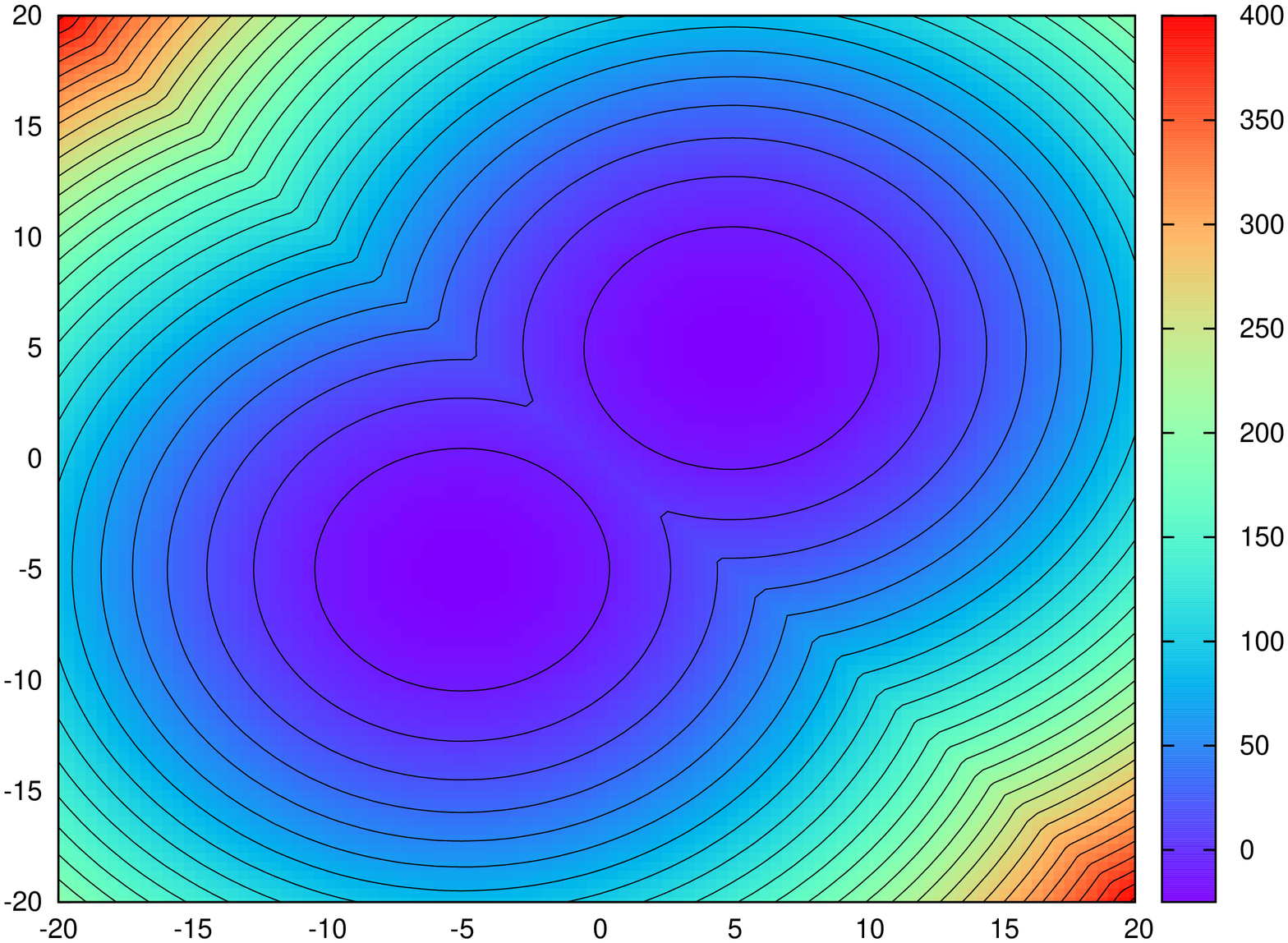}}
 \vspace{-.1cm}
  \centerline{\footnotesize{(a)}}\medskip
\end{minipage}
\hfill
 \begin{minipage}[b]{0.49\linewidth}
   \centering
   \centerline{\includegraphics[scale=0.35]{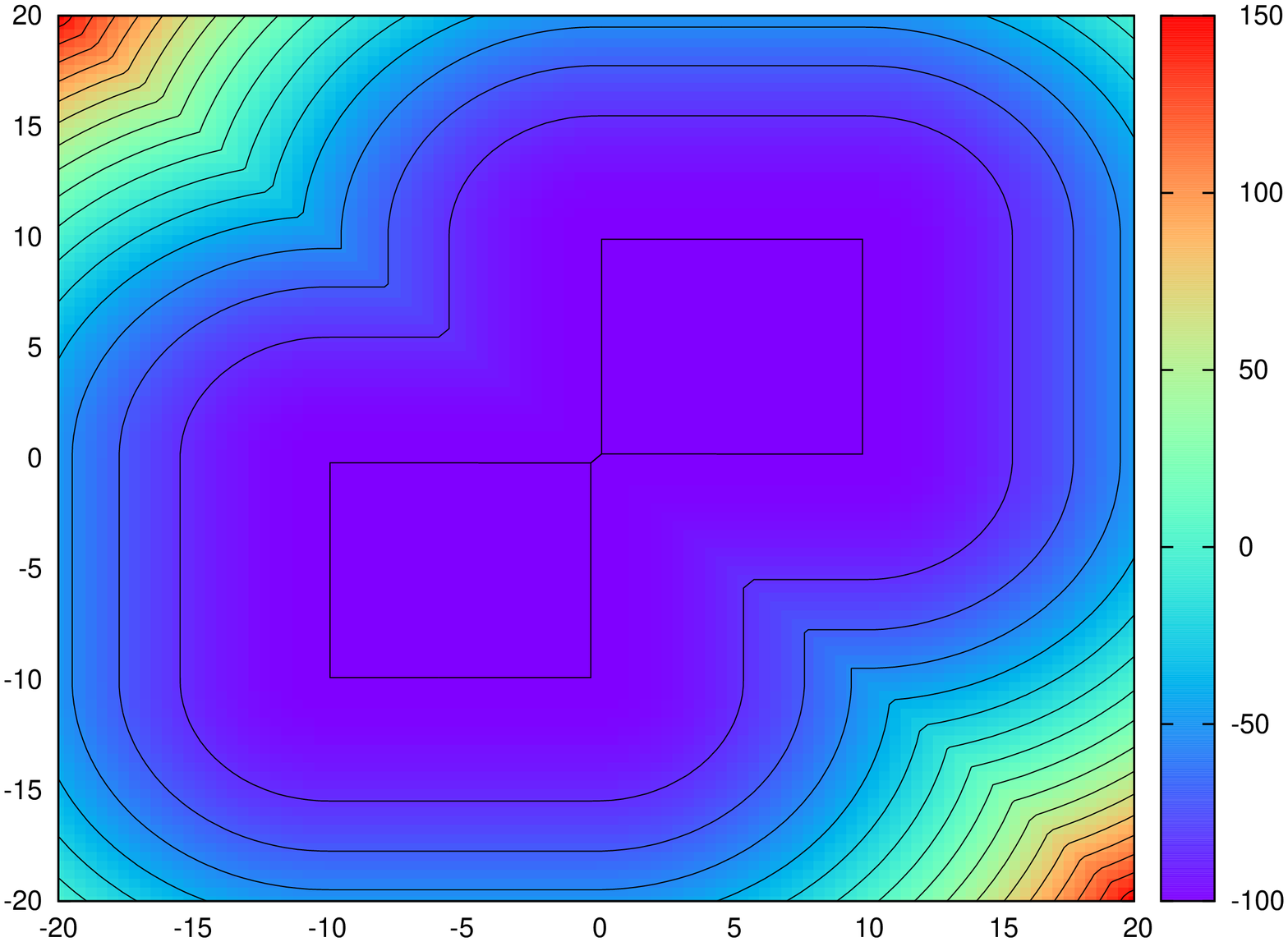}}
   \vspace{-.1cm}
   \centerline{\footnotesize{(b)}}\medskip
 \end{minipage}
\begin{minipage}[b]{0.49\linewidth}
  \centering
  \centerline{\includegraphics[scale=0.35]{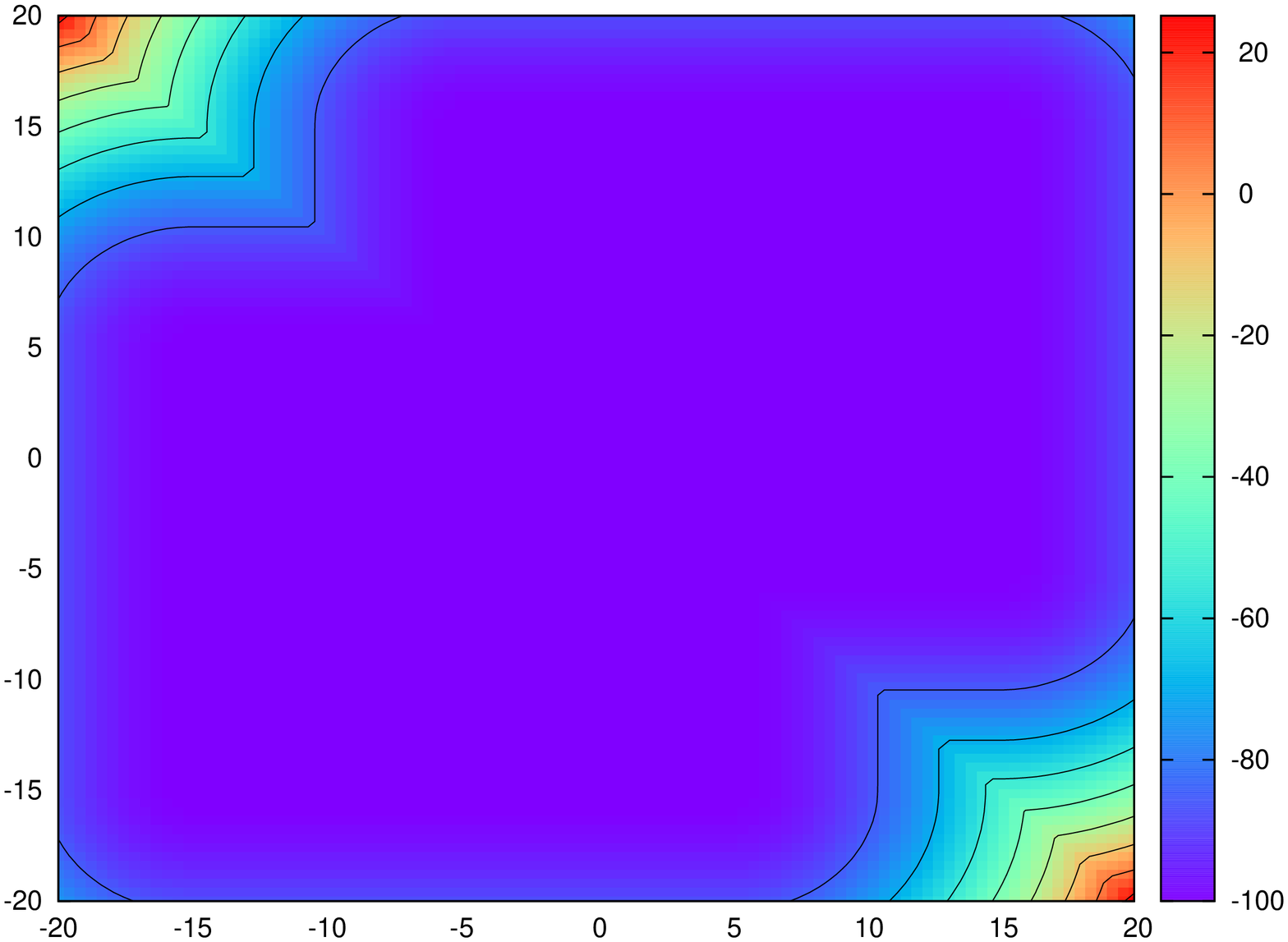}}
  \vspace{-.1cm}
  \centerline{\footnotesize{(c)}}\medskip
\end{minipage}
\begin{minipage}[b]{0.49\linewidth}
  \centering
  \centerline{\includegraphics[scale=0.35]{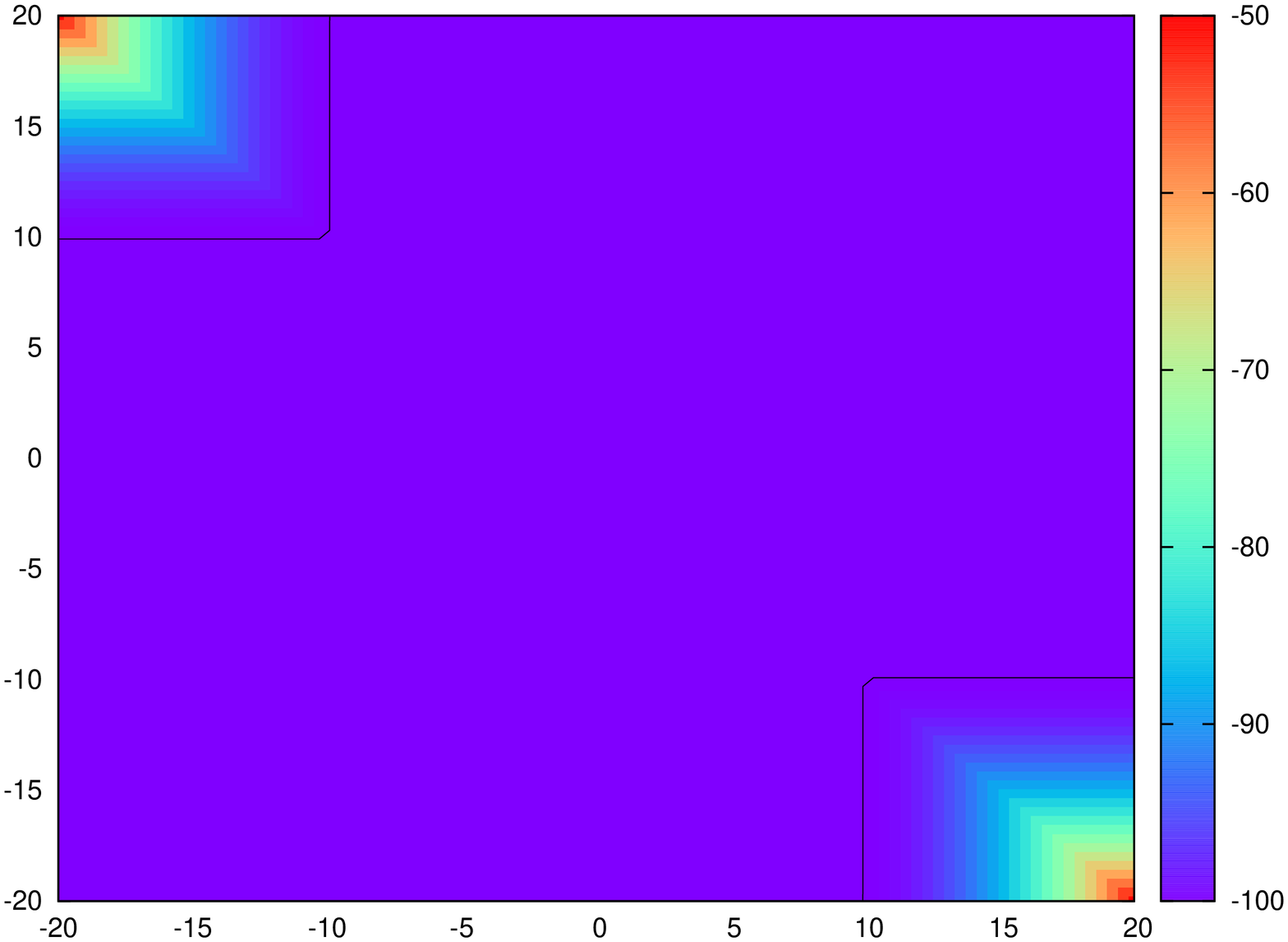}}
  \vspace{-.1cm}
  \centerline{\footnotesize{(d)}}\medskip
\end{minipage}

  \caption{
    Evaluation of the solution $\phi((x_1,x_2,0,0,0,0,0,0)^\dagger,t)$ of
    the HJ-PDE with initial data  $J=
    \min{}\left(\frac{1}{2}\|\cdot\|_2^2 - \langle b, \cdot\rangle,
    \frac{1}{2}\|\cdot\|_2^2 + \langle b, \cdot\rangle\right)$ with
  $b=(1,1,1,1,1,1,1,1)^\dagger$ and Hamiltonian $H=\|\cdot\|_1$, for
  $(x_1,x_2) \,\in\,     [-20,20]^2$ for   different times $t$.
    Plots for $t=0, 5,10,15$ are
    respectively depicted in (a), (b), (c) and (d). The level
    lines multiple of 15 are superimposed   on the plots.
}
\label{fig.min_initial_data}
\end{figure}

\begin{figure}[th]
\begin{minipage}[b]{.49\linewidth}
  \centering
 \centerline{\includegraphics[scale=0.35]{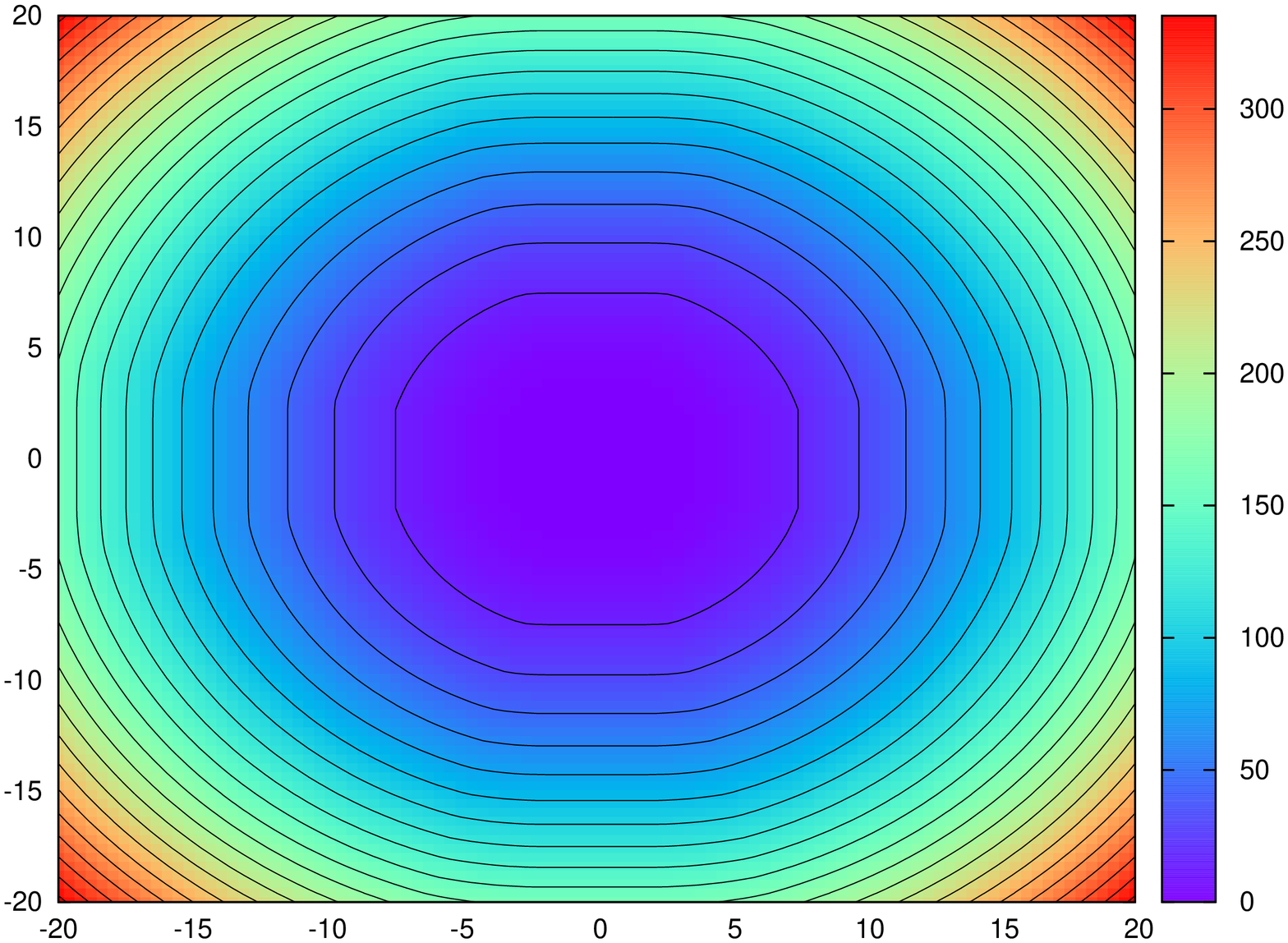}}
 \vspace{-.1cm}
  \centerline{\footnotesize{(a)}}\medskip
\end{minipage}
\hfill
 \begin{minipage}[b]{0.49\linewidth}
   \centering
 \centerline{\includegraphics[scale=0.35]{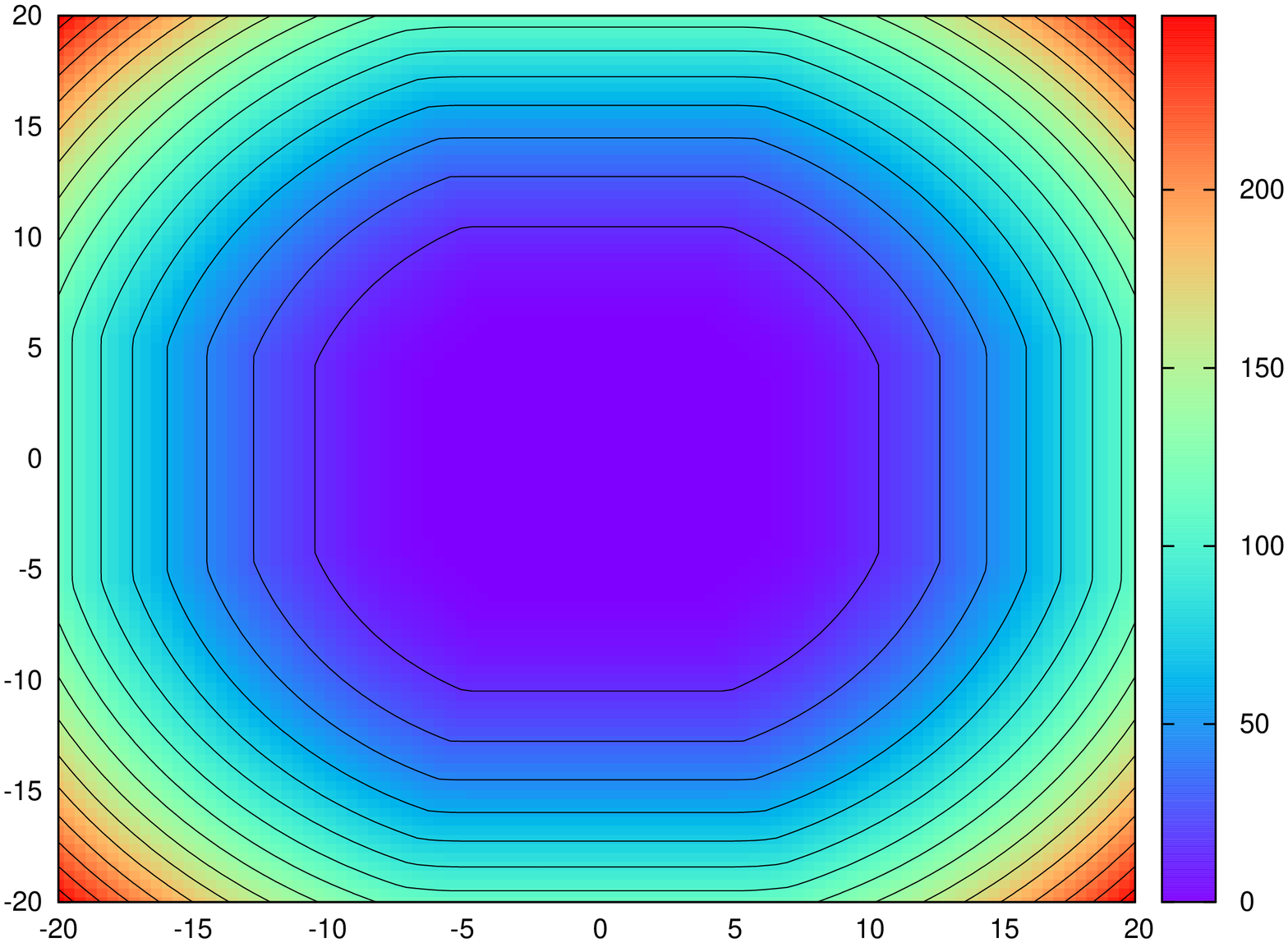}}
   \vspace{-.1cm}
   \centerline{\footnotesize{(b)}}\medskip
 \end{minipage}
\begin{minipage}[b]{0.49\linewidth}
  \centering
 \centerline{\includegraphics[scale=0.35]{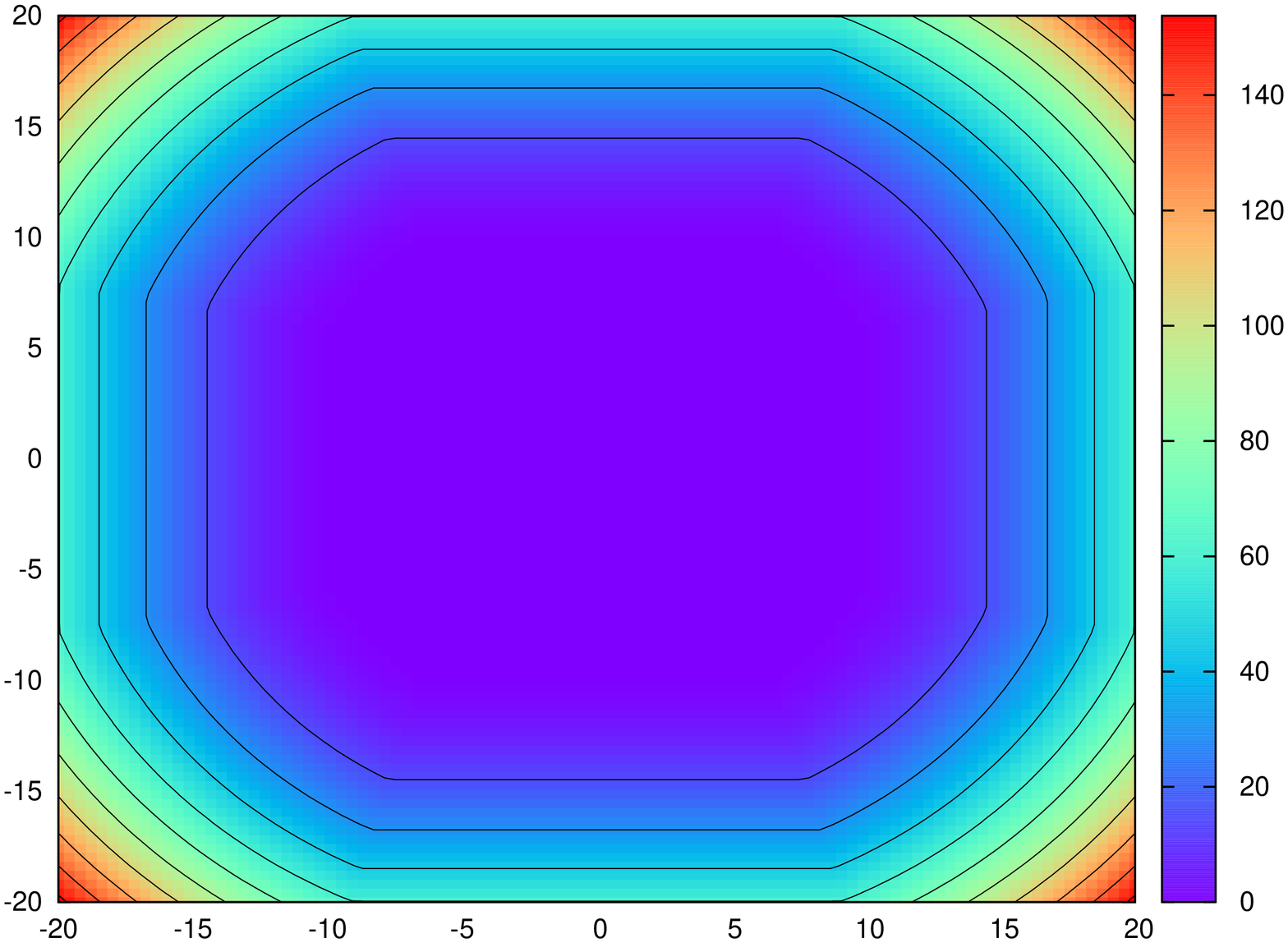}}
  \vspace{-.1cm}
  \centerline{\footnotesize{(c)}}\medskip
\end{minipage}
\begin{minipage}[b]{0.49\linewidth}
  \centering
 \centerline{\includegraphics[scale=0.35]{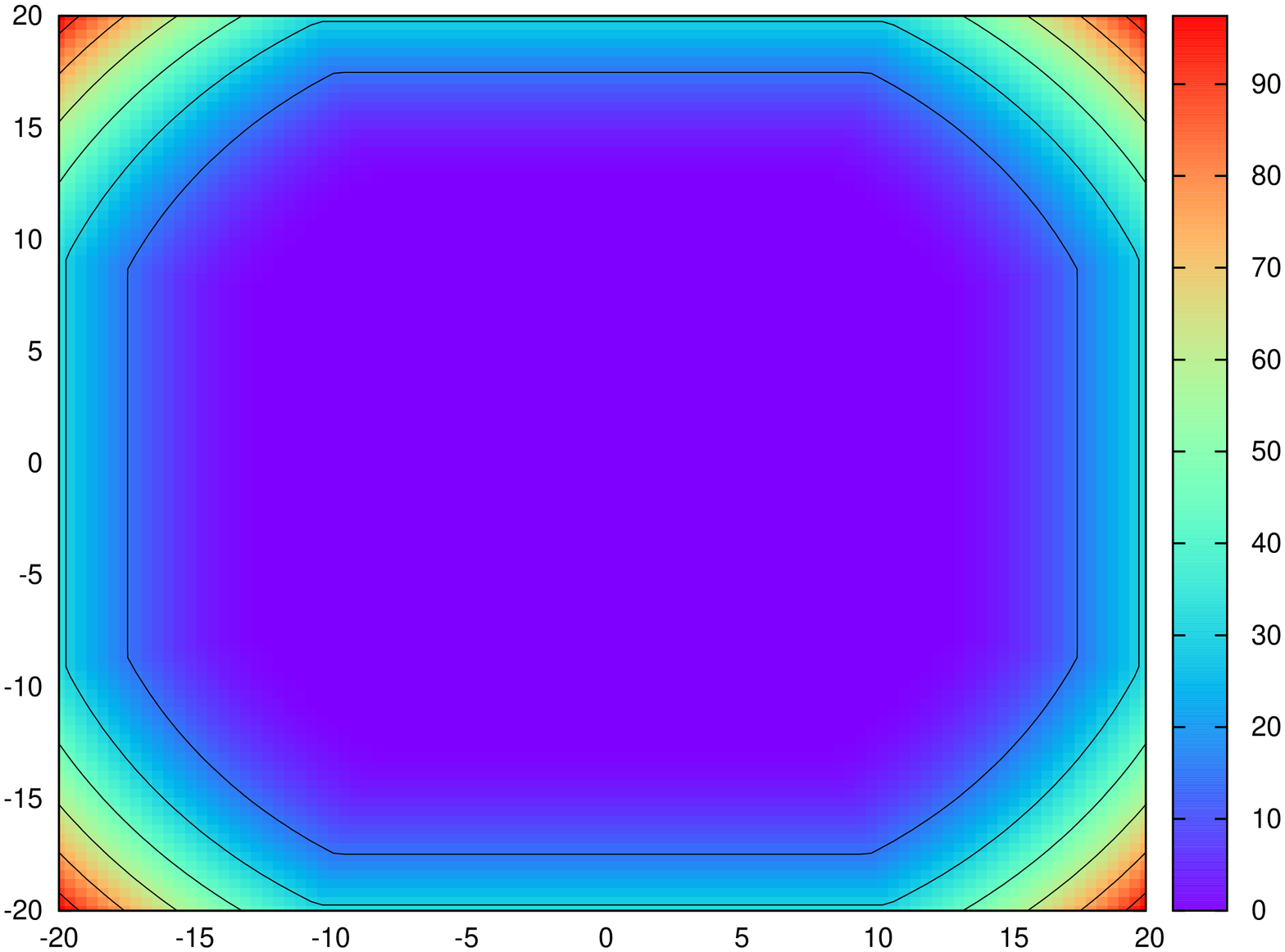}}
  \vspace{-.1cm}
  \centerline{\footnotesize{(d)}}\medskip
\end{minipage}

  \caption{
    Evaluation of the solution $\phi((x_1,x_2,0,0,0,0,0,0)^\dagger,t)$ of
    the HJ-PDE with initial data  $J=\frac{1}{2}\|\cdot\|_1^2$
    and Hamiltonian $H=\min{}\left(\|\cdot\|_1, \sqrt{\langle \cdot,
      \frac{4}{3}D \cdot\rangle} \right) $, for $(x_1,x_2) \,\in\,
    [-20,20]^2$ for   different times $t$.
    Plots for $t=2, 5,9,12$ are
    respectively depicted in (a), (b), (c) and (d). The level
    lines multiple of 15 are superimposed   on the plots.
}
\label{fig.min_hamiltonians}
\end{figure}

 We now consider solutions of HJ PDEs in dimension $n=8$
  on a 2-dimensional grid. We evaluate $\phi(x_1, x_2,0,0,0,0,0,0)$
 with $x_i \in \cup_{k=0,\dots,99}\{-20 + k    \frac{40}{99} \}$ for
 $i=1,2$.
Figures \ref{fig.J_linfty_2__H__l2},  \ref{fig.J_l1_2__H__l1} and
\ref{fig.J_l1_2__H__DIAG_NORM}  depict the solutions with initial data
$J=\frac{1}{2}\|\cdot\|_\infty^2$,
$J=\frac{1}{2}\|\cdot\|_1^2$ and $J=\frac{1}{2}\|\cdot\|_1^2$, and
with Hamiltonians $H= \|\cdot\|_2$, $H = \|\cdot\|_1$.
and $H=\sqrt{\langle \cdot, D \cdot\rangle}$ for various times,
respectively.
 Figure~\ref{fig.min_initial_data} and
  Figure~\ref{fig.min_hamiltonians} illustrate the max/min-plus
  algebra results described in
  Section~\ref{subsec.extension_future_work}. Figure~\ref{fig.min_initial_data} depicts the HJ solution for  the initial data $J =
J=    \min{}\left(\frac{1}{2}\|\cdot\|_2^2 - \langle b, \cdot\rangle,
    \frac{1}{2}\|\cdot\|_2^2 + \langle b, \cdot\rangle\right)$ with
  $b=(1,1,1,1,1,1,1,1)^\dagger$, and $H=\|\cdot\|_1$ for various
  times. Figure~\ref{fig.min_hamiltonians} depicts the HJ solution
  for various time with
  $J= \frac{1}\|\cdot\|_2^2$ and
$H=\min{}\left(\|\cdot\|_1, \sqrt{\langle \cdot,
      \frac{4}{3}D \cdot\rangle} \right) $.



\begin{table}[ht]
\begin{center}
  \begin{tabular}{r|c|c|c|c|c}
    n &  $\|y\|_1$ & $\|y\|_2$ & $\|y\|_\infty$  & $\|y\|_D$ & $\|y\|_A$\\
    \hline
    4  & 6.36e-08 & 1.20e-07 & 2.69e-07 & 7.00e-07 & 8.83e-07 \\
    8  & 6.98e-08 & 1.28e-07 & 4.89e-07 & 1.07e-06 & 1.57e-06 \\
    12 & 8.72e-08 & 1.56e-07 & 7.09e-07 & 1.59e-06 & 2.23e-06 \\
    16 & 9.24e-08 & 1.50e-07 & 9.92e-07 & 2.04e-06 & 2.95e-06 \\
  \end{tabular}
\end{center}
\caption{Time results in seconds for the average time per call for
  evaluting the solution of the HJ-PDE with the initial data $J =
  \frac{1}{2}\|\cdot\|_2^2$, several Hamiltonians and
various  dimensions~$n$.}
\label{table.l2_2}
\end{table}

\begin{table}[ht]
\begin{center}
  \begin{tabular}{r|c|c|c|c|c}
    n &  $\|y\|_1$ & $\|y\|_2$ & $\|y\|_\infty$  & $\|y\|_D$ & $\|y\|_A$\\
    \hline
    4   & 1.79e-06 & 1.53e-06 & 1.84e-06 & 4.88e-06 & 7.77e-06 \\
    8   & 3.77e-06 & 2.31e-06 & 3.50e-06 & 9.73e-06 & 1.92e-05 \\
    12  & 6.31e-06 & 3.14e-06 & 5.54e-06 & 1.44e-05 & 2.91e-05 \\
    16  & 9.61e-06 & 3.88e-06 & 8.22e-06 & 1.80e-05 & 4.04e-05 \\
  \end{tabular}
\end{center}
\caption{Time results in seconds for the average time per call for
  evaluting the solution of the HJ-PDE with the initial data $J =
  \frac{1}{2}\|\cdot\|_\infty^2$, several Hamiltonians and various
  dimensions~$n$.}
\label{table.linfty_2}
\end{table}

\begin{table}[ht]
\begin{center}
  \begin{tabular}{r|c|c|c|c|c}
    n &  $\|y\|_1$ & $\|y\|_2$ & $\|y\|_\infty$  & $\|y\|_D$ & $\|y\|_A$\\
    \hline
    4   & 2.86e-06 & 4.42e-06 & 9.17e-06 & 1.79e-05 & 1.97e-05 \\
    8   & 9.85e-06 & 1.63e-05 & 4.38e-05 & 9.37e-05 & 1.09e-04 \\
    12  & 2.35e-05 & 3.84-05  & 1.19e-04 & 2.63e-04 & 3.24e-04 \\
    16  & 4.35e-05 & 7.03e-05 & 2.46e-04 & 5.19e-04 & 6.92e-04 \\
  \end{tabular}
\end{center}
\caption{Time results in seconds for the average time per call for
  evaluting the solution of the HJ-PDE with the initial data $J =
  \frac{1}{2}\|\cdot\|_1^2$, several Hamiltonians and
  various dimensions $n$.}
\label{table.l1_2}
\end{table}

\begin{table}[ht]
\begin{center}
  \begin{tabular}{r|c|c|c|c|c}
    n &  $\|y\|_1$ & $\|y\|_2$ & $\|y\|_\infty$  & $\|y\|_D$ & $\|y\|_A$\\
    \hline
    4   & 3.62e-07 & 5.19e-07 & 9.35e-07 & 2.79e-06 & 3.50e-06 \\
    8   & 3.83e-07 & 5.25e-07 & 1.42e-06 & 4.40e-06 & 5.75e-06 \\
    12  & 4.97e-07 & 6.62e-07 & 1.73e-06 & 5.70e-06 & 7.98-06 \\
    16  & 5.92e-07 & 6.88e-07 & 2.27e-06 & 6.64e-06 & 1.04e-05 \\
  \end{tabular}
\end{center}
\caption{Time results in seconds for the average time per call for
  evaluting the solution of the HJ-PDE with the initial data $J =
  \frac{1}{2}\langle\cdot, D^{-1}\cdot \rangle$, several Hamiltonians and
  various dimensions $n$.}
\label{table.diag_2}
\end{table}

\begin{table}[h]
\begin{center}
  \begin{tabular}{r|c|c|c|c|c}
    n &  1 core &  4 cores & 8 cores  & 16 cores\\
    \hline
    4   & 1.11e-05 & 2.81e-06 & 1.56e-06  &  8.36e-07   \\
    8   & 4.77e-05 & 1.33e-05 & 6.81e-06  &  3.48e-06         \\
    12  & 1.35e-04 & 3.90e-05 & 1.94e-05  &  9.90e-06  \\
    16  & 3.24e-04 & 8.76e-05 & 4.40e-05  &  2.22e-05  \\
  \end{tabular}
\end{center}
\caption{Time results in seconds for the average time per call for
  evaluting the solution of the HJ-PDE with the initial data $J =
  \frac{1}{2}\|\cdot\|_1^2$, and the Hamiltonian
  $H=\|\cdot\|_\infty$, for various dimensions, and several cores.}
\label{table.parallel}
\end{table}



\newpage

\section{Conclusion}
\label{sec.conclusion}

We have designed algorithms which enable us to solve certain
Hamilton-Jacobi equations very rapidly. Our algorithms not only
evaluate the solution but also computes the gradient of the
solution. These include equations
arising in control theory leading to Hamiltonians which are convex and
positively homogeneous of degree $1$. We were motivated by ideas coming
from compressed sensing; we borrowed algorithms devised to solve
$\ell_1$ regularized problems which are known to rapidly converge.
We apparently extended this fast convergence to include convex
positively 1-homogeneous regularized problems.

There are no grids involved. Instead of complexity which is
exponential in the dimension of the problems, which is typical of grid
based methods, ours appears to be polynomial in the dimension with
very small constants. We can evaluate the solution on a laptop at
about $10^{-4} - 10^{-8}$ seconds per evaluation for fairly high
dimensions.  Our algorithm requires very low
  memory and is totally parallelizable which suggest that it is
  suitable for low energy embedded systems.
We have chosen to restrict the presentation of the numerical
experiments to norm-based Hamiltonians and we
  emphasize that our approach naturally extends to more elaborate
  positively 1-homogeneous Hamiltonians (using the min/max algebra
  results as we did for instance).

As an important step in this procedure we have also derived an equally
fast method to find a closest point
lying on $\Omega$, a finite union of compact convex sets $\Omega_i$,
such that $\Omega = \cup_{i}^{k}\Omega_i$ has a nonempty
interior, to a given point. 

We can also solve certain so called fast marching
\cite{tsitsiklis.95.itac} and fast sweeping \cite{TCOZ} problems equally
rapidly in high dimensions.  If we wish to
find $\psi:\Rn\to\R$ with, say $\psi=0$ on the boundary of a set
$\Omega$ defined above, satisfying
$$
\begin{cases}
\|\nabla_x \psi(x)\|_2 \;=\; 1 & \mbox{ in } \Rn,\\
\psi(x) \;=\; 0 & \mbox{ for any } x \in
(\Omega \setminus \interior{\Omega}),
\end{cases}
$$
then, we can solve for $u:\Rn\times
[0,+\infty) \to \R$
$$
\frac{\partial u}{\partial t} (x,t)
\,+\,  \|\nabla_x u(x)\|_2 \;=\; 0  \quad \mbox{ in } \Rn \times
(0,+\infty),
$$
with initial data
$$
\begin{cases}
u(x,0) < 0 \mbox{ for any } x \in \mbox{int } \Omega, \\
u(x,0) > 0 \mbox{ for any } x \in (\Rn \setminus \Omega),\\
u(x,0) = 0 \mbox{ for any } x \in (\Omega \setminus \mbox{int } \Omega),
\end{cases}
$$
and locate the zero level set of $u(\cdot, t) = 0$ for any given
$t>0$. Indeed any $x \in \{y \in \Rn \;|\; u(y,t) \;=\;0\}$
satisfies $\psi(x) = t$.

Of course the same approach could be used for any
 convex, positively 1-homogeneous Hamiltonian $H$ (instead of
 $\|\cdot\|_2$), e.g., $H \,=\,  \|\cdot\|_1$.
This will give us results related to computing the Manhattan distance.

We expect to extend our work as follows:
\begin{enumerate}
\item We will do experiments involving linear controls, allowing $x$
  and $t>0$ dependence while the Hamiltonian $(p,x,t) \mapsto H(p,x,t)$ is
  still convex and   positively 1-homogeneous  in $p$. The procedure
  was described in section \ref{sec.introduction_optimal_control}.
\item We will extend our fast computation of the projection in several
  ways. We will consider in detail the case of polyhedral regions
  defined by the intersection of sets $\Omega_i \,=\, \{x\in\Rn\;|\;
  \langle a_i, x \rangle - b_i \,\leq\,0\}$, $a_i,b_i \,\in\,\Rn$,
  $\|a_i\|_2=1$, for $i=1,\dots,k$. This is of interest in linear
  programming (LP) and related problems. We expect to develop
  alternate approaches to several issues arising in LP, including
  rapidly finding the existence and location of a feasible point.
\item We will consider nonconvex but positively 1-homogeneous
  Hamiltonians. These arise in: differential games as well as in the problem
  of finding a closest point on the boundary of a given compact convex
  set $\Omega$, to an arbitrary point in the {\it interior} of $\Omega$.
\item As an example of a nonconvex Hamiltonians we consider the
  following problems arising in differential games \cite{ES}. We need
  to solve the following scalar problem  for any $z\in\Rn$ and
any $\alpha >0$
\begin{displaymath}
\min_y \left\{\frac{1}{2} \|y-x\|_2^2 - \alpha\|y\|_1\right\}.
\end{displaymath}

It is easy to see that the minimizer is the stretch$_1$ operator which
we define for any $i=1,\dots,n$ as:
\begin{equation}
 \label{eq.stretch_def}
 \left(\hbox{stretch}_1 (x,\alpha)\right)_i \;=\;
\begin{cases}
 x_i+\alpha & \hbox{if } \ x_i > 0, \\
 0     & \hbox{if }  x_i=0,\\
 x_i-\alpha  & \hbox{if } \ x_i < 0.
\end{cases}
\end{equation}
We note that the discontinuity in the minimizer will lead to a jump
in the derivatives $(x,t)\mapsto \frac{\partial\varphi}{\partial
  x_i}(x,t)$, which is no surprise, given that this interface
associated with the equation
\begin{displaymath}
\frac{\partial\varphi}{\partial t}(x,t) - \sum_{i=1}^n \left|\frac{\partial
\varphi}{\partial x_i}(x,t)\right| \;=\; 0,
\end{displaymath}
and the previous initial data, will move inwards, and characteristics
will intersect.  The solution $\varphi(x,t)$ will remain locally Lipschitz
continuous, even though a point inside the ellipsoid may be equally
close to two points on the boundary of the original ellipsoid in the
Manhattan metric.
So we are solving
\begin{align}
\varphi(x,t) &= -\frac{1}{2} - \min_{v\in \Rn}\left\{\frac{1}{2}
\sum_{i=1}^n a_i^2 v_i^2 - t \sum_{i=1}^n |v_i| +
\langle x,v\rangle\right\} \nonumber \\
&= -\frac{1}{2} + \frac{1}{2} \sum_{i=1}^n \frac{x_i^2}{a_i^2} -
\min_{v\in \Rn}
\left\{\frac{1}{2} \sum_i^n a_i^2 \left(v_i -
\frac{x_i}{a^2}\right)^2 - t \sum_{i=1}^n |v_i|\right\} \nonumber \\
&= -\frac{1}{2} + \frac{1}{2} \sum_{i=1}^n \frac{(|x_i| + t)^2}{a_i^2} \nonumber
\end{align}
The zero level set disappears when $t \geq \max_{i} a_i$ as it should.

For completeness, we also consider the nonconvex optimization problem
\begin{displaymath}
\min_{v\in \Rn} \left\{\frac{1}{2} \|v-x\|_2^2 - \alpha\,\|v\|_2\right\}.
\end{displaymath}
Its minimizer is given by the $\mbox{stretch}_2$ operator formally defined by
$$
\mbox{stretch}_2(x, \alpha) \;=\;
\begin{cases}
  x + \alpha \frac{x}{\|x\|_2}  & \mbox{ if } x\neq 0,\\
  \alpha \theta  & \mbox{ with } \| \theta \|_2=1 \mbox{ if } x=0.
\end{cases}
$$

This formula, although multivalued at $x=0$, is useful to solve the
following problem: move the unit
sphere inwards with normal velocity $1$. The solution  comes from
finding the zero level set of
\begin{align}
\varphi(x,t) &= -\min_{v\in\Rn} \left\{\frac{|v|_2^2}{2} - t\|v\|_2 -
  \langle x,v\rangle\right\} -\frac{1}{2}
\nonumber \\
&= -\min_{v\in \Rn} \left\{\frac{1}{2} \|v-x\|_2^2 - t\|v\|_2\right\}
+\frac{1}{2} (\|x\|_2^2 - 1) \nonumber \\
&= -\frac{1}{2} t^2 + t\|x\|_2 \left(1 + \frac{t}{\|x\|_2}\right) +
\frac{1}{2} (\|x\|_2^2 - 1) \nonumber \\
&= \frac{1}{2} (\|x\|_2 + t)^2 - \frac{1}{2} \nonumber
\end{align}
and, of course, the zero level set is the set of $x$ satisfying $\|x\|_2
= t-1$ if $t \leq 1$ and the zero level set vanishes for $t > 1$.

\end{enumerate}


\appendix
\section{Gauge}
Now suppose furthermore that we wish $H$ to be nonnegative, i.e.,
$H(p) \geq 0$ for any $p \in \Rn$. Under this additional assumption,
we shall see that the Hamiltonian $H$ is not only the support function
of $C$ but also the gauge of convex set that we will characterize.
Let us define the gauge $\gamma_\Omega:\Rn \to \Rn \cup \{+\infty\}$ of
a closed convex set $D$ containing the origin (\cite[Def. 1.2.4,~p.
202]{hiriart-lemarechal.96.book-vol1})
\begin{equation}
  \nonumber
  \gamma_D(x)
 \;=\;
\inf \left\{ \lambda > 0 \;|\: x \in \lambda D\right\},
\end{equation}
where $\gamma_D(x) = + \infty$ if $ x \notin \lambda D$ for all
$\lambda > 0$. Using \cite[Thm. 1.2.5 (i), p
203]{hiriart-lemarechal.96.book-vol1}, $\gamma_D$ is lower
semicontinous, convex and positively 1-homogeneous. If $0 \in
\interior{D}$ then \cite[Thm. 1.2.5
(ii), p. 203]{hiriart-lemarechal.96.book-vol1} the gauge $\gamma_D$ is finite
everywhere, i.e., $\gamma_D:\Rn\to\R$, we recall that $\interior D$
denotes the interior of the set $D$. Thus, taking $H = \gamma_D$
satisfies our requirements. Note that if we further assume that $D$ is
symmetric (i.e., for any $d \in D$ then $-d \in D$) then $\gamma_D$ is
a seminorm. If in addition we wish  $H(p) > 0 $ for any $p \in \Rn
\setminus \{0\}$, then $D$ has to be a compact set
(\cite[Corollary. 1.2.6, p. 204]{hiriart-lemarechal.96.book-vol1}). For the
latter case, a symmetric $D$ implies that $\gamma_D$ is actually a
norm.

We now describe the connections between the sets $D$ and $C$, the
gauge $\gamma_D$, the support function of $C$ and nonnegative, convex,
positively 1-homogeneous Hamiltonians $H:\Rn\to \R$. Given a closed
convex set $\Omega$ of $\Rn$ containing the origin, we denote by
$\Omega^\circ$ the polar set of $\Omega$ defined by
\begin{equation}
\nonumber
\Omega^\circ \;=\;
\left\{
s \in \Rn \;|\; \langle s, x\rangle \leq 1 \mbox{ for all } x \in \Omega\}
\right\}.
\end{equation}
We have that $(\Omega^\circ)^\circ \,=\, \Omega$.
Using \cite[Prop. 3.2.4, p. 223]{hiriart-lemarechal.96.book-vol1} we
have that the gauge $\gamma_D$ is the support function of $D^\circ$.
In addition, \cite[Corollary 3.2.5,
p. 233]{hiriart-lemarechal.96.book-vol1} states that the support
function of $C$ is the gauge of $C^\circ$. We can take $D \,=\,
C^\circ$ and see that $H$ can be expressed as a gauge. Also, using
\cite[Thm. 1.2.5 (iii), p. 203]{hiriart-lemarechal.96.book-vol1} we
have that $C^\circ \,=\, \{ p \in \Rn \,|\, H(p) \leq 1\}$.


\section*{Acknowledgements} The authors deeply thank Gary Hewer
(Naval Air Weapon Center, China Lake) for fruitful discussions,
carefully reading drafts, and helping us to improve the paper.



\begin{thebibliography}{99}
\bibitem{akian.06.book_chapter}
M. Akian, R. Bapat, and S. Gaubert.
\newblock{\it Max-plus algebras.}
\newblock{In L. Hogben, editor, Handbook of Linear Algebra (Discrete
  Mathematics and Its Applications), volume 39. Chapman \& Hall/CRC,
  2006. Chapter 25.}

\bibitem{aubin.84.book}
J.-P. Aubin and A.~Cellina.
\newblock {\it Differential Inclusions}.
\newblock Springer-Verlag, Berlin, 1984.

\bibitem{RB1} R. Bellman, {\it Adaptive Control Processes, a Guided
    Tour}, Princeton U. Press (1961).

\bibitem{RB2} R. Bellman, {\it Dynamic Programming}, Princeton
  U. Press, (1957).

\bibitem{brezis.73.book}
H.~Brezis.
\newblock {\it Op\'erateurs maximaux monotones et semi-groupes de
  contractions dans les espaces de Hilbert}.
\newblock North-Holland, 1973.

\bibitem{CRT} E.J. Candes, J. Romberg, and T. Tao, {\it Exact Signal Reconstruction from Highly Incomplete Frequency Information}, IEEE Trans. on Information Theory, v.52, (2), (2008), pp. 489-509.

 \bibitem{chambolle.09.ijcv}
 A.~Chambolle and J.~Darbon.
 \newblock{\it On total variation minimization and surface evolution using
   parametric maximum flows}.
 \newblock { International Journal of Computer Vision}, 84(3):288--307, 2009.

\bibitem{chambolle.11.jmiv}
A.~Chambolle and T. Pock.
\newblock{\it A First-Order Primal-Dual Algorithm for Convex Problems
  with Applications to Imaging}.
\newblock{Journal of Mathematical Imaging and Vision}, 41(1):120-145, 2011.

\bibitem{cheng-tsai.08.jcp}
L.T. Cheng and Y.H. Tsai.
\newblock{\it Redistancing by flow of time dependent eikonal equation.}
\newblock{ J. Comput. Phys. 227(8), 2008}

\bibitem{combettes.07.siam-opt}
P.L. Combettes and J.-C. Pesquet.
\newblock {\it Proximal thresholding algorithm for minimization over orthonormal  bases.}
\newblock { SIAM Journal on Optimization}, 18(4):1351--1376, November 2007.



\bibitem{CEL} M.G. Crandall, L.C. Evans, and P.-L. Lions, {\it Some Properties of Viscosity Solutions of Hamilton-Jacobi Equations}, Trans. AMS, v.282, (2), (1984), pp. 487-502.

\bibitem{CL} M.G. Crandall and P.-L. Lions, {\it Viscosity Solutions
    of Hamilton-Jacobi Equations}, Trans. AMS, v.277, (1), (1983), pp.
  1-42.

\bibitem{darbon.13.cam}
J.~Darbon.
\newblock{\it On Convex Finite-Dimensional Variational Methods in Imaging Sciences, and Hamilton-Jacobi Equations.}
\newblock { SIAM Journal on Imaging Sciences 8:4, 2268-2293, 2015}.



\bibitem{daubechies.04.cpam}
I.~Daubechies, M.~Defrise, and C.~De~Mol.
\newblock{\it  An iterative thresholding algorithm for linear inverse
  problems with  a sparsity constraint.}
\newblock { Comm. Pure and Appl. Math.},
  57(11):1413--1457, 2004.


\bibitem{Di}
E.W. Dijkstra, {\it A Note on Two Problems in Connexion
  `    with Graphs}, Numer. Math., v.1, (1959), pp. 269-271.

\bibitem{Dol} I.C. Dolcetta, {\it Representation of Solutions of
Hamilton-Jacobi Equations}, Nonlinear Equations: Methods, Models and
Applications, (2003), pp. 79-90.

\bibitem{Do} D.L. Donoho, {\it Compressed Sensing}, IEEE Trans. on Image Information Theory, v.52, (4), (2006), pp. 1289-1305.

\bibitem{ekeland.76.book}
I.~Ekeland and R.~Temam.
\newblock {\it Convex Analysis and Variational Problems}.
\newblock North-Holland, Amsterdam, 1976.

\bibitem{Ev} L.C. Evans, {\it Partial Differential Equations},
  Grad. Studies in Mathematics, v. 19, AMS (2010).

\bibitem{ES} L.C. Evans and P.E. Souganidis, {\it Differential Games
    and Representation Formulas for Solutions of Hamilton-Jacobi
    Isaacs Equations}, Indiana U. Math. J., v.38, (1984), pp. 773-797.


\bibitem{figueiredo.98.asolimar}
M.A.T. Figueiredo and R.D. Nowak.
\newblock{\it Bayesian wavelet-based signal estimation using
  non-informative   priors.}
\newblock In { Signals, Systems amp; Computers, 1998. Conference Record of
  the Thirty-Second Asilomar Conference on}, volume~2, IEEE
Piscataway, NJ 1998, pp. 1368--1373.

\bibitem{fleming.97.pisa}
W. H. Fleming
\newblock{\it Deterministic nonlinear filtering.}
\newblock { Annali della Scuola Normale Superiore di Pisa - Classe di
  Scienze}, 25(3-4):435-454, 1997.


\bibitem{GM} R. Glowinski and A. Marrocco, {\it Sur l'approximation
    par \'el\'ements finis d'ordre, et la resolution par
    p\'enalzation-dualit\'e, d'une classe de problem\'emes de
    Dirichlet non-lin\'eares}, C.R. hebd.
S\'eanc. Acad. Sci., Paris 278, s\'erie A, (1974), pp. 1649-1652.

\bibitem{GO} T. Goldstein and S. Osher, {\it The Split Bregman Method
    for $L_1$ Regularized Problems}, SIAM J. on Imaging Science, v.2,
  (2), 2009, pp. 323-343.


\bibitem{He} M.R. Hestenes, {\it Multiplier and Gradient Methods}, J.
  of Optimization Theory and Applications, v.4, (3), (1969), pp.
  303-320.

\bibitem{hiriart.98.book}
J.-B. Hiriart-Urruty.
\newblock {\it Optimisation et Analyse Convexe}.
\newblock Presse Universitaire de France, 1998.

\bibitem{hiriart-lemarechal.96.book-vol1}
J.-B. Hiriart-Urruty and C.~{Lemar{\'e}chal}.
\newblock {\it Convex Analysis and Minimization Algorithms Part I}.
\newblock Springer Verlag, Heidelberg, 1996.

\bibitem{hiriart-lemarechal.96.book-vol2}
J.-B. Hiriart-Urruty and C.~{Lemar{\'e}chal}.
\newblock {\it Convex Analysis and Minimization Algorithms Part II}.
\newblock Springer Verlag, Heidelberg, 1996.

\bibitem{eckstein-bertsekas.92.mpa}
J. Eckstein and D.P. Bertsekas
\newblock{\it On the Douglas-Rachford splittingmethod and the proximal
  point algorithm for maximal monotone operators}
\newblock{ Math. Program., Ser. A 55(3), 293\u2013318 (1992)}.

\bibitem{Ho} E. Hopf, {\it Generalized Solutions of Nonlinear
    Equations of the First Order}, J. Math. Mech., v.14, (1965),
  pp. 951-973.

\bibitem{HS} C. Hu and C.-W. Shu, {\it A Discontinuous Galerkin Finite
    Element Method for Hamilton-Jacobi Equations}, SIAM
  J. Sci. Comput., v.21, (2), (1999), pp. 666-690.

\bibitem{kurzhanski-varaiya.14.book}
A. B. Kurzhanski and P. Varaiya
\newblock{\it Dynamics and Control of Trajectory Tubes: Theory and
  Computation.}
\newblock{Birkh\"auser, 2014 edition.}

\bibitem{lions.79.sna}
P.-L. Lions and B.~Mercier.
\newblock{\it Splitting algorithms for the sum of two nonlinear operators.}
\newblock { SIAM Journal on Numerical Analysis}, 16(6):964--979, 1979.

\bibitem{lions-rochet.86.ams}
P.-L. Lions and J-C Rochet.
\newblock{\it Hopf formula and multitime Hamilton-Jacobi equations.}
\newblock{ Proc. of American Mathmatical Socity, 96, no.1 (1986),
  79--84.}

\bibitem{mceneaney.06.book}
W. M. McEneaney
\newblock{\it Max-Plus Methods for Nonlinear Control and Estimation .}
\newblock{Birkh\"auser, 2006 edition.}

\bibitem{MBT} I.M. Mitchell, A.M. Bayen and C.J. Tomlin, {\it A Time
    Dependent Hamilton-Jacobi Formulation of Reachable Sets for
    Continuous Dynamic Games}, IEEE Trans on Automatic Control, v.50,
  171, (2005), pp. 947-957.

\bibitem{MT} I.M. Mitchell and C.J. Tomlin, {\it Overapproximating
    Reachable Sets by Hamilton-Jacobi Projections}, J. Sci. Comput.,
  v.19, (1-3), (2003), pp. 323-346.

\bibitem{moreau.65.proximite}
J.-J. Moreau.
\newblock{\it Proximit\'e et dualit\'e dans un espace hilbertien.}
\newblock { Bulletin Spc. Math. France}, 93, pp. 273--299, 1965.

\bibitem{OOTT} A. Oberman and S. Osher and R. Takei and R. Tsai,
\newblock{\it Numerical methods for anisotropic mean curvature flow based on a
  discrete time variational formulation,}
\newblock{Commun. Math. Sci., v.9, no. 3, pp.637-662, 2011.}

\bibitem{OM} S. Osher and B. Merriman, {\it The Wulff Shape as the Asymptotic Limit of a Growing Crystalline Interface}, Asian J. Mathematics, v.1, (3), (1997), pp. 560-571.

\bibitem{OSe} S. Osher and J.A. Sethian, {\it Fronts Propagating with
    Curvature Dependent Speech: Algorithms Based on Hamilton-Jacobi
    Formulations}, J. Comput. Phys., v.79, (1), (1988), pp. 12-49.

\bibitem{OS} S. Osher and C.-W. Shu, {\it High Order Essentially
    Nonosicllatory Schemes for Hamilton-Jacobi Equations}, SIAM
  J. Numer. Analysis, v.28, (4), (1991), pp. 907-922.

\bibitem{OY} S. Osher and W. Yin, {\it Error Forgetting of Bregman Iteration}, J. Sci. Comput., v.54, (2), (2013), pp. 684-695.

\bibitem{rockafellar.70.book}
R.T. Rockafellar.
\newblock{\it  Convex Analysis}.
\newblock{ Princeton Landmarks in
Mathematics. Princeton University Press, Princeton (1997).
Reprint of the 1970 original, Princeton Paperbacks.}

\bibitem{rockafellar.76.siam-co}
R.T. Rockafellar
\newblock{ \it Monotone operators and the proximal point algorithm.}
\newblock{SIAM J. Control Optim. 14(5), 877-898 (1976).}


\bibitem{rockafellar.84.book}
R.T. Rockafellar.
\newblock{\it Network Flows and Monotropic Optimization}.
\newblock  Reprint of the 1984 original. Athena Scientific, 1998.

\bibitem{TCOZ}
Y.-H. R. Tsai, L.-T. Cheng, S. Osher, and H.-K. Zhao.
{\it Fast Sweeping Algorithms for a Class of Hamilton--Jacobi Equations
}
SIAM Journal on Numerical Analysis 41:2, 673-694, 2003.


\bibitem{teboulle.97.siam-op}
M. Teboulle
\newblock{Convergence of Proximal-like Algorithms.}
\newblock{SIAM J. Optimization 7 (1997), 1069-1083.}


\bibitem{tsitsiklis.95.itac}
J. N. Tsitsiklis
\newblock{Efficient algorithms for globally optimal trajectories.}
\newblock{Automatic Control, IEEE Transactions on , vol.40, no.9,
  pp.1528-1538, Sep 1995}


\bibitem{winkler.03.book}
G.~Winkler.
\newblock { \it{Image Analysis,
      Random Fields and Dynamic Monte Carlo Methods}}.
\newblock  Applications of mathematics. 2nd edition, Springer-Verlag, 2006.



\bibitem{YOGD} W. Yin, S. Osher, D. Goldfarb, and J. Darbon, {\it
    Bregman Iterative Algorithms for $\ell_1$ Minimization
with Applications to Compressed Sensing}, SIAM J. on Imaging Sci., v.1
(1), (2008), pp. 143-168.

\bibitem{ZS} Y.T. Zhang and C.-W. Shu, {\it High Order WENO Schemes for Hamilton-Jacobi Equations on Triangular Meshes}, SIAM J. on Sci. Comput., v.24, (3), (2013), pp. 1005-1030.

\bibitem{zhu-chan.08.cam}
M. Zhu and T.F. Chan,
{\it An Efficient Primal-Dual Hybrid Gradient Algorithm For Total
  Variation Image Restoration},
UCLA CAM Report 08-34, 2008.

\end{thebibliography}
\end{document}